\documentclass{article}
\usepackage{amsfonts}
\begin{document}
\title {\bf  Hypoelliptic Dunkl equations
in the space of distributions on $\mathbb{R}^d$}
\date{}
\author{\bf Khalifa TRIM\`ECHE\\
{\small Faculty of Sciences of Tunis, Department of
Mathematics},\\ {\small CAMPUS,  1060 Tunis, Tunisia}\\ email :
{\small khlifa.trimeche@fst.rnu.tn}} \maketitle
\begin{abstract}
In this paper we define and study the Dunkl convolution product
and the   Dunkl transform on spaces of distributions on $
\mathbb{R}^d$. By using the main results obtained, we study the
hypoelliptic Dunkl convolution equations in the space of
distributions.
\end{abstract}
\section{Introduction} We consider the differential-difference
operators $T_j, j = 1, 2,\cdots,d$, on $\mathbb{R}^d$ introduced
by C.F.Dunkl in [5]. These operators are very important  in pure
mathematics and in Physics. They provide a useful tool in the
study of special functions with root systems [6,9,4].  Moreover
the commutative algebra generated by these operators has been used
in the study of certain exactly solvable models of quantum
mechanics, namely the Calogero-Sutherland-Moser models, which deal
with systems of identical particles in a one dimensional space
(see [10,15,16]).

C.F.Dunkl has proved in [7] that there exists a unique isomorphism
$V_k$ from the space of homogeneous polynomials $\mathcal{P}_n$ on
$\mathbb{R}^d$ of degree $n$ onto itself satisfying the
transmutation relations $$T_jV_k = V_k \frac{\partial}{\partial
x_j},\quad j = 1, 2,\cdots,d,\eqno{(1.1)}$$ and $$V_k(1) =
1.\eqno{(1.2)}$$ This operator is called Dunkl intertwining
operator. It has been extended to an isomorphism from
$\mathcal{E}(\mathbb{R}^d)$ (the space of $C^\infty$-functions on
$\mathbb{R}^d)$ onto itself satisfying the relations (1.1) and
(1.2) (see[23]).

The operator $V_k$ possesses the integral  representation
$$\forall\; x \in \mathbb{R}^d,\quad V_k(f)(x) =
\int_{\mathbb{R}^d}f(y)d\mu_x(y), \quad f \in
\mathcal{E}(\mathbb{R}^d),\eqno{(1.3)}$$ where $\mu_x$ is a
probability measure on $\mathbb{R}^d$ with \mbox{support} in the
closed ball $B(0, \|x\|)$ of center $0$ and radius $\|x\|$ (see
[19, 23]).

We have shown in [23] that for each $x \in \mathbb{R}^d$, there
exists a unique distribution $\eta_x$ in
$\mathcal{E}'(\mathbb{R}^d)$ (the space of distributions on
$\mathbb{R}^d$ of compact \mbox{support}) with \mbox{support} in
$B(0, \|x\|)$ such that $$V_k^{-1}(f)(x) = \langle
\eta_x,f\rangle, f \in \mathcal{E}(\mathbb{R}^d).$$

We have studied also in [23] the transposed operator ${}^tV_k$ of
the operator $V_k$. It has the integral representation $$\forall\;
y \in \mathbb{R}^d,\;\; {}^tV_k(f)(y) =
\int_{\mathbb{R}^d}f(x)d\nu_y(x).\eqno{(1.4)} $$ where $\nu_y$ is
a positive measure on $\mathbb{R}^d$ with \mbox{support} in the
set $\{x \in \mathbb{R}^d /\;\|x\| \geq \|y\|\}$ and $f$ in
$\mathcal{D}(\mathbb{R}^d)$ (the space of $C^\infty$-functions on
$\mathbb{R}^d$ with compact \mbox{support}).\\ This operator is
called Dual Dunkl intertwining operator.

We have proved in [23] that the operator ${}^tV_k$ is an
isomorphism from $\mathcal{D}(\mathbb{R}^d)$ onto itself,
satisfying the transmutation relations $$\forall\; y \in
\mathbb{R}^d,\;\;\; {}^tV_k(T_jf)(y) = \frac{\partial}{\partial
y_j} {}^tV_k(f)(y),\quad j = 1, 2,\cdots,d.\eqno{(1.5)}$$

Using the operator $V_k$, C.F.Dunkl has defined in [7] the Dunkl
kernel $K$ by $$\forall\; x \in \mathbb{R}^d,\; \forall\; z \in
\mathbb{C}^d,\; K(x, -iz) = V_k(e^{-i\langle
.,z\rangle})(x).\eqno{(1.6)}$$ Using this kernel C.F.Dunkl has
introduced in [7] a Fourier transform $\mathcal{F}_D$ called Dunkl
transform.\\ \hspace*{5mm}By using the operators $V_k$ and
$\,^{t}V_k$ we have defined in [25] the Dunkl translation
operators and we have determined their properties. With the aid of
these operators we define and study in this work the Dunkl
convolution product on spaces of distributions. We present also
the properties of the Dunkl transform of distributions. The
results obtained have permitted to characterize the hypoelliptic
Dunkl convolution equations in the space of distributions in terms
of their Dunkl transform. This characterization was first given by
L.Ehrenpreis [8] and next by L.H\"ormander [11] in the case of the
classical Fourier transform on $ \mathbb{R}^d$. In [1][2] the
authors have studied this characterization for the Hankel, Jacobi
and Ch\'ebli-Trim\`eche transforms. We remark that their proof of
the existence of a parametrix is complicated. In this paper we
give a very simple proof of this result for the Dunkl transform on
$ \mathbb{R}^d$, which can also be applied to the cases of the
preceding transforms.

\section{The eigenfunction of the Dunkl operators}

In this section we collect some notations and results on Dunkl
operators and the Dunkl kernel (see [6, 7, 12, 13,14]).
\subsection{Reflection Groups, Root Systems and Multiplicity
Functions}

We consider $\mathbb{R}^d$ with the euclidean scalar product
$\langle .,.\rangle$ and $\|x\| = \sqrt{\langle x , x\rangle}$. On
$\mathbb{C}^d, \|.\|$ denotes also the standard Hermitian norm,
while $\langle z, w\rangle = \sum^d_{j=1}z_j \overline{w_j}$ .

For $\alpha \in \mathbb{R}^d \backslash \{0\}$, let
$\sigma_\alpha$ be the reflection in the hyperplan $H_\alpha
\subset \mathbb{R}^d$ orthogonal to $\alpha$, i.e.
$$\sigma_\alpha(x) = x - \left(\frac{2\langle
\alpha,x\rangle}{\|\alpha\|^2}\right)\alpha .\eqno{(2.1)}$$ A
finite set $R \subset \mathbb{R}^d \backslash \{0\}$ is called a
root  system if $R \cap \mathbb{R}\alpha = \{\pm \alpha\}$ and
$\sigma_\alpha R = R$ for all $\alpha \in R$. For a given  root
system $R$ the reflections $\alpha_\alpha, \alpha \in R$, generate
a finite group $W \subset O(d)$, the reflection group associated
with $R$. All reflections in $W$ correspond to suitable pairs of
roots. For a given $\beta \in \mathbb{R}^d \backslash \cup_{\alpha
\in R }H_\alpha$, we fix the positive subsystem $R_+ = \{\alpha
\in R ; \langle \alpha, \beta\rangle > 0\}$, then for each $\alpha
\in R$ either $\alpha \in R_+$ or $-\alpha \in R_+$.

A function $k : R \rightarrow \mathbb{C}$ on a root system $R$ is
called a multiplicity function if it is invariant under the action
of the associated reflection group $W$. If one regards $k$ as a
function on the corresponding reflections, this means that $k$ is
constant on the conjugacy classes of reflections in $W$. For
abbreviation, we introduce the index $$\gamma = \gamma(R) =
\sum_{\alpha \in R_+} = \sum_{\alpha \in R_+}
k(\alpha).\eqno{(2.2)}$$ Moreover, let $\omega_k$ denotes the
weight function $$\omega_k(x) = \prod_{\alpha \in R_+}|\langle
\alpha,x\rangle|^{2k(\alpha)}.\eqno{(2.3)}$$ which is
$W$-invariant and homogeneous of degree $2\gamma$.

For $d = 1$ and $W= \mathbb{Z}_2$, the multiplicity function $k$
is a single parameter denoted  also $k $ and $$\forall\; x \in
\mathbb{R},\;\; \omega_k(x) = |x|^{2k}.\eqno{(2.4)} $$ We
introduce the Mehta-type constant $$c_k = \left(
\int_{\mathbb{R}^d}
e^{-\|x\|^2}\omega_k(x)dx\right)^{-1}.\eqno{(2.5)}$$ which is
known for all Coxeter groups $W$ (see [5, 9])
\subsection{Dunkl Operators and Dunkl kernel}

The Dunkl operators $T_j, j = 1,\cdots,d$, on $\mathbb{R}^d$,
associated with the finite reflection group $W$ and the
multiplicity function $k$, are given for a function $f$ of class
$C^1$ on $\mathbb{R}^d$ by $$T_jf(x) = \frac{\partial}{\partial
x_j}f(x) + \sum_{x \in R_+} k(\alpha)\alpha_j \frac{f(x) -
f(\sigma_\alpha(x))}{\langle\alpha,x\rangle}.\eqno{(2.6)} $$ In
the case $k = 0$, the $T_j, j = 1, 2,\cdots,d$, reduce to the
corresponding partial derivatives. In this paper, we will assume
throughout that $k \geq 0$ and $\gamma > 0$.

For $f$ of class $C^1$ on $\mathbb{R}^d$ with compact
\mbox{support} and $g$ of class $C^1$ on $\mathbb{R}^d$ we have
$$\int_{\mathbb{R}^d}T_jf(x) g(x) \omega_k(x)dx =
-\int_{\mathbb{R}^d} f(x) T_j g(x) \omega_k(x) dx,\quad j = 1,
2,\cdots, d.\eqno{(2.7)}$$ For $y \in \mathbb{R}^d$, the system
$$\left\{ \begin{array}{ll} T_j u(x, y) &= y_j u(x, y),\quad j =
1, 2,\cdots,d,\\ u(0,y) &= 1,
\end{array}\right.\eqno{(2.8)}$$
admits a unique analytic solution on $\mathbb{R}^d$, denoted by
$K(x,y)$ and called Dunkl kernel.

This kernel has a unique holomorphic extension to $\mathbb{C}^d
\times \mathbb{C}^d$.\\ {\bf Example 2.1.}

If $d = 1$ and $W = \mathbb{Z}_2$, the Dunkl kernel is given by
$$K(z,t) = j_{\gamma-1/2} (izt) + \frac{zt}{2\gamma+1}
j_{\gamma+1/2} (izt),\quad z, t \in \mathbb{C},\eqno{(2.9)}$$
where for $\alpha `\geq - 1/2, j_\alpha$ is the normalized Bessel
function defined by $$j_\alpha(u) = 2^\alpha \Gamma(\alpha + 1)
\frac{J_\alpha(u)}{u^\alpha} = \Gamma(\alpha +1) \sum^\infty_{n=0}
\frac{(-1)^n(u/2)^{2n}}{n!\Gamma(n+\alpha +1)}, \quad u \in
\mathbb{C},\eqno{(2.10)} $$ with $J_\alpha$  the Bessel function
of first kind and index $\alpha$ (see [7]).

The Dunkl kernel possesses the following properties.
\begin{itemize}
\item[(i)] For $z, t \in \mathbb{C}^d$, we have $K(z,t) = K(t,z), K(z,0) =
1$, and $K(\lambda z,t) = K(z, \lambda t)$ for all $\lambda \in
\mathbb{C}$.
\item[(ii)] For all $\nu \in \mathbb{Z}^d_+, x \in \mathbb{R}^d$,
and $z \in \mathbb{C}^d$ we have $$|D^\nu_z K(x,z)| \leq
\|x\|^{|\nu|} \exp \left[\max_{w \in W} \langle wx, Re
z\rangle\right].\eqno{(2.11)}$$ In particular $$|D^\nu_z K(x,z)|
\leq \|x\|^{|\nu|} \exp [\|x\| \|Re z\|]],\eqno{(2.12)} $$
$$|K(x,z)| \leq \exp [\|x\| \|Re z\|],\eqno{(2.13)}$$ and for all
$x, y \in \mathbb{R}^d$ : $$|K(ix,y)| \leq 1, \eqno{(2.14)}$$ with
$$D^\nu_z = \frac{\partial^{|\nu|}} {\partial z_1^{\nu_1} \cdots
\partial z^{\nu_d}_d } \mbox{ and } |\nu| = \nu_1 + \cdots+
\cdots + \nu_d.$$ \item[(iii)] For all $x, y \in \mathbb{R}^d$
 and $w \in W$ we have
 $$K(-ix, y) = \overline{K(ix, y)} \mbox{ and } K(wx, wy) =
  K(x,y).\eqno{(2.15)}$$
 \item[(iv)] The function $K(x,z)$ admits for all $x \in
 \mathbb{R}^d$ and $z \in \mathbb{C}^d$ the following Laplace type
 integral representation
 $$K(x,z) = \int_{\mathbb{R}^d} e^{\langle y,z\rangle}d\mu_x(y),
 \eqno{(2.16)}$$
 where $\mu_x$ is a probability measure on $\mathbb{R}^d$ with
 \mbox{support} in the closed ball $B(0, \|x\|)$ of center $0$ and radius
 $\|x\|$. and we have
 $$\mbox{supp}\, \mu_x \cap \{y \in \mathbb{R}^d / \|y\| = \|x\|\} \neq
 \emptyset.\eqno{(2.17)}$$ More precisely the measure $\mu_x$ satisfies \\
 \hspace*{5mm} - $\mbox{supp}\, \, \mu_x$ is contained in $co \{wx, \; w \in
 W\}$ the convex hull of the orbit of $x$ under $W$.\\  \hspace*{5mm} -
 $\mbox{supp}\, \mu_x \cap \{wx, \; w \in
 W\} \neq
 \emptyset.$
\end{itemize}
(see [19]).\\
 {\bf Remark 2.1}

When $d = 1$ and $W = \mathbb{Z}_2$, the relation (2.16) is of the
form $$K(x,z) = \frac{\Gamma(\gamma+1/2)}
{\sqrt{\pi}\Gamma(\gamma)}|x|^{-2\gamma}\int^{|x|}_{-|x|}(|x| -
y)^{\gamma-1}(|x| + y)^\gamma e^{yz}dy.\eqno{(2.18)}$$ Then in
this case the measure $\mu_x$ is given for all $x \in
\mathbb{R}\backslash \{0\}$ by: $$d\mu_x(y) = {\cal K}(x,y) dy,$$
with $${\cal K}(x,y) =
\frac{\Gamma(\gamma+1/2)}{\sqrt{\pi}\Gamma(\gamma)}|x|^{-2\gamma}(|x|
-y)^{\gamma-1}(|x|+y)^\gamma 1_{]-|x|, |x|[}(y),\eqno{(2.19)} $$
where $1_{]-|x|, |x|[}$ is the characteristic function of the
interval $]-|x|, |x|[$.

We remark that by change of variables,  the relation (2.18) takes
the following form $$\forall\; x \in \mathbb{R}^d, \; \forall\; z
\in \mathbb{C}^d,\;\; K(x,z) = \frac{\Gamma(\gamma + 1/2)}
{\sqrt{\pi}\Gamma(\gamma)}\int^1_{-1}e^{txz}(1-t^2)^{\gamma-1}(1+t)dt,\eqno{(2.20)}
$$
\section{The Dunkl intertwining operator and its dual}
{\bf Notation} We denote by

- $C(\mathbb{R}^d)$ (resp. $C_c(\mathbb{R}^d))$ the space of
continuous functions on $\mathbb{R}^d$ (resp. with compact
\mbox{support}).

- $C^p(\mathbb{R}^d)$ (resp. $C^p_c(\mathbb{R}^d))$ the space of
functions of class $C^p$ on $\mathbb{R}^d$ (resp. with compact
\mbox{support}).\\ We provide the preceding spaces with the
classical topology.\\ \hspace*{5mm}- $\mathcal{E}(\mathbb{R}^d)$
the space of $C^\infty$-functions on $\mathbb{R}^d$ equipped with
the topology of uniform convergence on all compact for the
functions and their derivatives.\\ \hspace*{5mm}-
$\mathcal{D}(\mathbb{R}^d)$ the space of $C^\infty$-functions on
$\mathbb{R}^d$ with compact \mbox{support}. We have $$
\mathcal{D}( \mathbb{R}^d) = \bigcup_{ a \geq 0}\mathcal{D}_a(
\mathbb{R}^d), $$ where $\mathcal{D}_a( \mathbb{R}^d)$ is the
space of $C^{\infty}$-functions on $ \mathbb{R}^d$, with support
in the closed ball $B(o,a)$ of center $o$ and radius $a$. \\ The
topology on $\mathcal{D}_a( \mathbb{R}^d)$ is defined by the
seminorms $$p_n(\psi) =
\renewcommand{\arraystretch}{0.5}
\begin{array}[t]{c}
\sup\\ {\scriptstyle |\mu| \leq n}\\ {\scriptstyle x \in B(o,a)}
\end{array}\renewcommand{\arraystretch}{1} |D^\mu \psi(x)|, \;  n
\in \mathbb{N},$$ where $$D^\mu = \frac{\partial^{|\mu|}}{\partial
x_1^{\mu_1}...\partial x_d^{\mu_d} }, \; \mu = (\mu_1,...,\mu_d)
\in \mathbb{N}^d.$$ These seminorms are equivalent to the
seminorms $$q_m(\psi) =
\renewcommand{\arraystretch}{0.5}
\begin{array}[t]{c}
\sup\\ {\scriptstyle |\mu| \leq m}\\ {\scriptstyle x \in B(o,a)}
\end{array}\renewcommand{\arraystretch}{1} |T^\mu \psi(x)|, \;  m
\in \mathbb{N},$$where $$T^\mu = T^{\mu_1}_1\;o\;
T_2^{\mu_2}o\cdots \;o T_d^{\mu_d}.$$ The space $ \mathcal{D}(
\mathbb{R}^d)$ equipped with the inductive limit topology is a
Fr\'echet space.\\ \hspace*{5mm} - $\mathcal{S}(\mathbb{R}^d)$ the
space of $C^\infty$-functions on $\mathbb{R}^d$ which are rapidly
decreasing as their derivatives. The topology on this space is
defined by  the seminorms $$P_{r,s}(\psi)
=
\renewcommand{\arraystretch}{0.5}
\begin{array}[t]{c}
\sup\\ {\scriptstyle |\mu| \leq r}\\ {\scriptstyle x \in
\mathbb{R}^d} \end{array}\renewcommand{\arraystretch}{1}(1 +
\|x\|^2)^s |D^\mu \psi(x)|, r, s \in \mathbb{N}.$$ These seminorms
are equivalent to the seminorms
 $$Q_{k,\ell}(\psi) =
\renewcommand{\arraystretch}{0.5}
\begin{array}[t]{c}
\sup\\ {\scriptstyle |\mu| \leq k}\\ {\scriptstyle x \in
\mathbb{R}^d} \end{array}\renewcommand{\arraystretch}{1}(1 +
\|x\|^2)^\ell |T^\mu \psi(x)|, k, \ell \in \mathbb{N}.$$ Equipped
with this topology $\mathcal{S}(\mathbb{R}^d)$ is a Fr\'echet
space.\\
 We consider also the following spaces.

- $\mathcal{E}'(\mathbb{R}^d)$ the space of distributions on
$\mathbb{R}^d$ with compact \mbox{support}. It is the topological
dual of $\mathcal{E}(\mathbb{R}^d)$.

- $\mathcal{S}'(\mathbb{R}^d)$ the space of tempered distributions
on $\mathbb{R}^d$. It is the topological dual of
$\mathcal{S}(\mathbb{R}^d)$.\vspace{5mm}

The Dunkl intertwining operator $V_k$ is defined on
$\mathcal{C}(\mathbb{R}^d)$ by $$\forall\; x \in \mathbb{R}^d,
V_k(f)(x) = \int_{\mathbb{R}^d}f(y)d\mu_x(y),\eqno{(3.1)}$$ where
$\mu_x$ is the measure given by the relation (2.16) (see [23]).

We have $$\forall\; x \in \mathbb{R}^d,\;\; \forall\; z \in
\mathbb{C}^d,\;\; K(x,z) = V_k(e^{<.,z>})(x).\eqno{(3.2)}$$ The
operator ${}^tV_k$ satisfying for $f$ in $C_c(\mathbb{R}^d)$ and
$g$ in $C(\mathbb{R}^d)$, the relation
$$\int_{\mathbb{R}^d}{}^tV_k(f)(y)g(y)dy =
\int_{\mathbb{R}^d}V_k(g)(x) f(x) \omega_k(x)dx.\eqno{(3.3)}$$ is
given by $$\forall\; y \in \mathbb{R}^d, {}^tV_k(f)(y) =
\int_{\mathbb{R}^d}f(x) d\nu_y(x),\eqno{(3.4)}$$ where $\nu_y$ is
a positive measure on $\mathbb{R}^d$ whose \mbox{support}
satisfies $$\mbox{supp}\, \nu_y \subset \{x \in \mathbb{R}^d /
\|x\|\geq \|y\|\} \mbox{ and } \mbox{supp}\, \nu_y \cap \{x \in
\mathbb{R}^d / \|x\| = \|y\|\} \neq \emptyset. \eqno{(3.5)} $$
This operator is called the dual Dunkl intertwining operator (see
[23]). \vspace{3mm}

 The following theorems give some properties
of the operators $V_k$ and ${}^tV_k$ (see [23]).\\ {\bf Theorem
3.1}
\begin{itemize}
\item[(i)] The operator $V_k$ is a topological isomorphism from
$\mathcal{E}(\mathbb{R}^d)$ onto itself satisfying the
transmutation relations $$\forall\; x \in \mathbb{R}^d,\;\;
T_jV_k(f)(x) = V_k\left(\frac{\partial}{\partial
y_j}f\right)(x),\;\; j = 1, 2,\cdots,d,\;\; f \in
\mathcal{E}(\mathbb{R}^d).\eqno{(3.6)}$$
\item[(ii)] For each $x \in \mathbb{R}^d$, there exists a unique
distribution $\eta_x$ in $\mathcal{E}'(\mathbb{R}^d)$ with
\mbox{support} in the closed ball $B(0, \|x\|)$ such that for all
$f$ in $\mathcal{E}(\mathbb{R}^d)$ we have $$V^{-1}_k(f)(x) =
\langle \eta_x, f\rangle.\eqno{(3.7)}$$
\end{itemize}
Moreover $$\mbox{supp}\, \eta_x \cap \{y \in \mathbb{R}^d / \|y\|
= \|x\|\} \neq \emptyset. \eqno{(3.8)} $$ {\bf Theorem 3.2}
\begin{itemize}
\item[(i)] The operator ${}^tV_k$ is a topological isomorphism
form $\mathcal{D}(\mathbb{R}^d)$ (resp.
$\mathcal{S}(\mathbb{R}^d))$ onto itself, satisfying the
transmutation relations $$\forall\; y \in \mathbb{R}^d,\;
{}^tV(T_jf)(y) = \frac{\partial}{\partial y_j}{}^tV(f)(y),\; j =
1, 2,\cdots,d, f \in \mathcal{D}(\mathbb{R}^d).\eqno{(3.9)} $$
\item[ii)] For all $f$ in $ \mathcal{D}( \mathbb{R}^d)$ we have
$$\mbox{supp} \, f \subset B(o,a) \Longleftrightarrow \mbox{supp}
\, \,^{t}V_k(f) \subset B(o,a).
 \eqno{(3.10)}$$
 where $ B(o,a)$ is the closed ball of center $o$ and radius $a > 0$.
\item[iii)] For each $y \in \mathbb{R}^d$, there exists a unique
distribution $Z_y$ in $\mathcal{S}'(\mathbb{R}^d)$ with
\mbox{support} in the set $\{x \in \mathbb{R}^d / \|x\| \geq
\|y\|\}$ such that for all $f$ in $\mathcal{D}(\mathbb{R}^d)$ we
have $${}^tV_k^{-1}(f)(y) = \langle Z_y, f\rangle.\eqno{(3.11)}$$
Moreover $$\mbox{supp}\, Z_y \cap \{ x \in \mathbb{R} / \|x\| =
\|y\|\} \neq \emptyset.\eqno{(3.12)}$$
\end{itemize}
{\bf Example 3.1}

When $d = 1$ and $W = \mathbb{Z}_2$, the Dunkl intertwining
operator $V_k$ is defined by (3.1) with for all $x \in \mathbb{R}
\backslash \{0\}, d\mu_x(y) = {\cal K}(x,y) dy$, where ${\cal K}$
given by the relation (2.19).

The dual Dunkl intertwining operator ${}^tV_k$ is defined by (3.4)
with $d\nu_y(x) = {\cal K}(x,y) \omega_k(x)dx$, where ${\cal K}$
and $\omega_k$ given respectively by the relations (2.19) and
(2.3).\\ {\bf Example 3.2}

The Dunkl intertwining operator $V_k$ of index $\gamma =
\sum^d_{i=1} \alpha_i, \alpha_i > 0$, associated with the
reflection group $\mathbb{Z}_2 \times \mathbb{Z}_2 \times \cdots
\times \mathbb{Z}_2$ on $\mathbb{R}^d$, is given for all $f$ in
$\mathcal{E}(\mathbb{R}^d)$ and for all $x \in \mathbb{R}^d$ by
$$V_k(f)(x) = \prod^d_{i=1}\left(\frac{\Gamma(\alpha_i + 1/2)}
{\sqrt{\pi}\Gamma(\alpha_i)}\right) \int_{[-1,1]d}f(t_1x_1,
t_2x_2,\cdots,t_dx_d)$$ \hspace{4cm}$ \times
\displaystyle{\prod^d_{i=1}}(1-t^2_i)^{\alpha_i-1}(1+t_i)dt_1\cdots
dt_d,$\hfill(3.13)\\ (see [27]).\vspace*{5mm}\\ {\bf{Definition
3.1.}} The dual Dunkl intertwining operator on ${\cal E'}(
\mathbb{R}^d)$ denoted also by $\,^{t}V_k$ is defined by $$\langle
\,^{t}V_k(S),\varphi\rangle = \langle S,V_k(\varphi)\rangle, \;
\varphi \in  \mathcal{E}( \mathbb{R}^d).\eqno{(3.14)} $$
\hspace*{5mm} The operator $\,^{t}V_k$ possesses the following
properties (See [25] p.26-27).\\ \hspace*{5mm} i) It is a
topological isomorphism from ${\cal E'}( \mathbb{R}^d)$ onto
itself. Its inverse is given by $$\langle
\,^{t}V_k^{-1}(S),\varphi\rangle = \langle
S,V_k^{-1}(\varphi)\rangle, \; \varphi \in  \mathcal{E}(
\mathbb{R}^d).\eqno{(3.15)} $$ \hspace*{5mm} ii) Let
$T_{f\omega_k}$ be the distribution of ${\cal E'}( \mathbb{R}^d)$
given by the function $f\omega_k$, with $f \in D( \mathbb{R}^d)$.
Then we have $$ \,^{t}V_k(T_{f\omega_k}) =
T_{\,^{t}V_k(f)}.\eqno{(3.16)} $$
 \hspace*{5mm} iii) Let
$T_{g}$ be the distribution of ${\cal E'}( \mathbb{R}^d)$ given by
the function $g$ in $D( \mathbb{R}^d)$. Then we have $$
\,^{t}V_k^{-1}(T_{g}) =
T_{\,^{t}V_k^{-1}(g)\omega_k}.\eqno{(3.17)} $$
\section{Dunkl transform}

In this section we define the Dunkl transform and we give the main
results satisfied by this transform (see [7, 13, 14]).\\ {\bf
Notations} We denote by

- $L^p_k(\mathbb{R}^d), p \in [1, + \infty]$, the space of
measurable functions on $\mathbb{R}^d$ such that
\begin{eqnarray*}
\|f\|_{k,p} &=& \left(\int_{\mathbb{R}^d}|f(x)|^p
\omega_k(x)dx\right)^{1/p}  < + \infty,\;\; \mbox{ if } 1\leq p <
+ \infty,\\ \|f\|_{k,\infty} &=& \displaystyle{ess\sup_{x \in
\mathbb{R}^d}}|f(x)| < + \infty.
\end{eqnarray*}

- $H(\mathbb{C}^d)$ the space of entire functions on
$\mathbb{C}^d$ which are rapidly decreasing and of exponential
type. We have $$H(\mathbb{C}^d) = \bigcup_{a \geq 0}
H_a(\mathbb{C}^d) $$ where $H_a(\mathbb{C}^d)$ is the space of
entire functions $\Psi$ on $ \mathbb{C}^d$ satisfying $$\forall \,
m \in \mathbb{N}, \; \sup_{z \in \mathbb{C}^d}(1+||z||^2)^m
|\Psi(z)|e^{-a||Im z||} < +\infty. $$ \vspace{3mm}\\
 The Dunkl transform of a function $f$ in
$\mathcal{D}(\mathbb{R}^d)$ is given by $$\forall\; y \in
\mathbb{R}^d,\;\; \mathcal{F}_D(f)(y) =
\int_{\mathbb{R}^d}f(x)K(x, -iy)\omega_k(x)dx.\eqno{(4.1)}$$ This
transform has the following properties.
\begin{itemize}
\item[i)] For $f$ in $L^1_k(\mathbb{R}^d)$ the function
$\mathcal{F}_D(f)$ belongs to $C(\mathbb{R}^d)$,
 tends to zero as $t$ goes to infinity, and
 we have $\|\mathcal{F}_D(f)\|_{k, \infty}
\leq \|f\|_{k,1}$.
\item[(ii)] Let $f$ be in $\mathcal{D}(\mathbb{R}^d)$. If $ \check{f}(x) = f(-x)$
and $f_w(x) = f(wx)$ for $x \in \mathbb{R}^d$,  $w \in W$, then
for all $y \in \mathbb{R}^d$ we have $$\mathcal{F}_D(
\check{f})(y) = \overline{\mathcal{F}_D(f)(y)} \mbox{ and }
\mathcal{F}_D(f_w)(y) = \mathcal{F}_D(f)(wy).\eqno{(4.2)}$$
\item[iii)] There is a one-to-one correspondence between
the space of all radial functions $f$ in
$L_{k}^{1}(\mathbb{R}^{d})$ and the space of integrable functions
$F$ on $[0,+\infty[$ with respect to the measure $\frac{r^{2
\gamma + d - 1} dr}{\Gamma(\gamma + \frac{d}{2})2^{ \gamma +
\frac{d}{2}}}$, via $$f(x) = F(||x||) = F(r), \; with \; r =
||x||.$$ Moreover, the Dunkl transform ${\cal F}_{D}(f)$ of $f$ is
related to the
 Fourier-Bessel transform ${\cal F}_{B}^{\gamma + \frac{d}{2} -
1}(F)$ of $F$ by $$
 \forall y \in \mathbb{R}^{d}, {\cal F}_{D}(f)(y) =
\frac{ 2^{\gamma + \frac{d}{2}}} {c_{k}} {\cal F}_{B}^{\gamma +
\frac{d}{2} - 1}(F)(||y||).\eqno{(4.3)}$$ The transform ${\cal
F}_{B}^{\gamma + \frac{d}{2} - 1}$ is  given by $$ \forall \lambda
\geq 0, \; {\cal F}_{B}^{\gamma + \frac{d}{2} - 1}(\lambda) =
\displaystyle
 \displaystyle \int_{0}^{\infty} g(r) j_{\gamma  + \frac{d}{2} - 1 }
(\lambda r) \frac{r^{  2 \gamma  + d - 1}}{\Gamma(\gamma +
\frac{d}{2})2^{ \gamma + \frac{d}{2}}}\; dr,\eqno{(4.4)}$$ with
$j_{\gamma  + \frac{d}{2} - 1 } (\lambda r) $  the normalized
Bessel function.  (See [18] p.585-589, and [24]).\\
\item[i$\nu$)] For all $f$ in $\mathcal{S}(\mathbb{R}^d)$ we have
$$\mathcal{F}_D(f) = \mathcal{F}\,o\,{}^tV_k(f),\eqno{(4.5)}$$
\end{itemize}
where $\mathcal{F}$ is the classical Fourier transform on
$\mathbb{R}^d$ given by $$\forall\; y \in \mathbb{R}^d,\;\;
\mathcal{F}(f)(y) = \int_{\mathbb{R}^d}f(x) e^{-i\langle
x,y\rangle}dx,\;\; f\in \mathcal{D}(\mathbb{R}^d),\eqno{(4.6)}$$
The following theorems are proved in [13, 14].\\ {\bf Theorem
4.1.} The transform $\mathcal{F}_D$ is a topological isomorphism
\begin{itemize}
\item[i)] from $\mathcal{D}(\mathbb{R}^d)$ onto
$\mathbb{H}(\mathbb{C}^d)$,
\item[ii)] from $\mathcal{S}(\mathbb{R}^d)$ onto itself.
\end{itemize}
The inverse transform is given by
 $$\forall\; x \in \mathbb{R}^d,\;\; \mathcal{F}^{-1}_D(h)(x)
=
\frac{c^2_k}{2^{2\gamma+d}}\int_{\mathbb{R}^d}h(y)K(x,iy)
\omega_k(y)dy.\eqno{(4.7)}$$ \noindent{\bf Remark 4.1}

Another proof of Theorem 4.1 is given in [25].\\ {\bf Theorem
4.2.} Let $f$ be in $L^1_k(\mathbb{R}^d)$ such that the function
$\mathcal{F}_D(f)$ belongs to $L^1_k(\mathbb{R}^d)$. Then we have
the following inversion formula for the transform $\mathcal{F}_D$
: $$f(x) = \frac{c^2_k}{2^{2\gamma+d}}\int_{\mathbb{R}^d}
\mathcal{F}_D(f)(y)K(x, iy)\omega_k(y)dy,\;\; a.e.\eqno{(4.8)}$$
{\bf Theorem 4.3.}
\begin{itemize}
\item[i)] \underline{Plancherel formula for}  $\mathcal{F}_D$.\\ For all $f$ in
$\mathcal{D}(\mathbb{R}^d)$ we have $$\int_{\mathbb{R}^d}|f(x)|^2
\omega_x(x)dx =
\frac{c^2_k}{2^{2\gamma+d}}\int_{\mathbb{R}^d}|\mathcal{F}_D(f)(y)|^2
\omega_k(y)dy.\eqno{(4.9)} $$
\item[ii)] \underline{Plancherel Theorem for} $\mathcal{F}_D$.\\ The
renormalized Dunkl transform $f\rightarrow 2^{-\gamma-d/2}  c_k
\mathcal{F}_D(f)$ can be uniquely extended to an isometric
isomorphism on $L^2_k(\mathbb{R}^d)$.
\end{itemize}
\section{Dunkl convolution product and Dunkl transform of distributions}
\subsection{Dunkl translation operators and Dunkl convolution product of functions}

The definitions and properties of Dunkl translation operators and
Dunkl convolution product of functions presented in this
subsection are given in the seventh section of [25] p. 33 - 37.

The Dunkl translation operators $\tau_x, x \in \mathbb{R}^d$, are
defined on $\mathcal{E}(\mathbb{R}^d)$ by $$\forall\; y \in
\mathbb{R}^d, \tau_x f(y) = (V_k)_x(V_k)_y[V^{-1}_k
(f)(x+y)].\eqno{(5.1)}$$ For $f$ in $\mathcal{D}(\mathbb{R}^d)$
the function $\tau_x f$ can  be expressed by using the dual Dunkl
intertwining operator as follows $$\forall\; y \in \mathbb{R}^d,
\tau_x f(y) = (V_k)_x ({}^tV^{-1}_k)_y
[{}^tV_k(f)(x+y)].\eqno{(5.2)}$$ Using the relations (5.1) and
(5.2) we deduce that the Dunkl translation operators  can also be
written in the following forms $$\forall\; (x,y) \in
\mathbb{R}^d\times\mathbb{R}^d, \tau_x f(y) = (V_k)\otimes
V_k)(\Psi) (x,y),\eqno{(5.3)}$$ with $$\Psi(x,y) =
V^{-1}_k(f)(x+y), \; f \in \mathcal{E}( \mathbb{R}^d),$$ and
$$\forall\; y \in \mathbb{R}^d, \tau_x f(y) = \,^{t}V_k^{-1}[
 \check{\mu}_x * \,^{t}V_k(f)](y),\eqno{(5.4)}$$ where $\mu_x$ is the
measure given by (2.16), $  \check{\mu}_x$ the measure defined by
$$\int_{ \mathbb{R}^d} g(y)d \check{\mu}_x(y) =  \int_{
\mathbb{R}^d} g(-y)d\mu_x(y), \; g \in C(
\mathbb{R}^d),\eqno{(5.5)}$$ and $ \check{\mu}_x * g$ the function
given by $$\forall\; y \in \mathbb{R}^d, \;  \check{\mu}_x * g(y)
= \int_{ \mathbb{R}^d} g(y-t)d \check{ \mu}_x(t).\eqno{(5.6)}$$
\hspace*{5mm} The operators $\tau_x$, $x \in \mathbb{R}^d$,
satisfy the properties\\ \hspace*{7mm}i) For all $ x \in
\mathbb{R}^{d} $, the operators $\tau_{x} $, is continuous from
$\mathcal{E}( \mathbb{R}^{d} )$ into itself.\\ \hspace*{5mm} ii)
The function $x \mapsto \tau_{x} $, is of class $C^{\infty}$ on
$\mathbb{R}^{d}.$\\ \hspace*{5mm} iii)
 For all $x, y \in \mathbb{R}^{d}$ and $ z \in {\mathbb{C}}^{d}$ we
have the product formula $$
 \tau_{x} K(y,z) = K(x,z) K(y,z).\eqno{(5.7)}
$$ \hspace*{6mm}$i\tau$) For all  $f$ in ${\cal
E}(\mathbb{R}^{d})$, we have $$ \tau_{x} f(0) =  f(x), \;
 \tau_{x} f(y) =  \tau_{y} f(x).\eqno{(5.8)}
$$ and $$ T_{j}(\tau_{x} f) = \tau_{x}(T_{j} f), \quad j = 1,...,
d.\eqno{(5.9)} $$ $$ (T_{j})_{x}(\tau_{x} f) = \tau_{x}(T_{j} f),
\quad j = 1,..., d.\eqno{(5.10)} $$ \hspace*{5mm} $\nu)$ For $f$
in $\mathcal{D}(\mathbb{R}^d)$ and $x \in \mathbb{R}^d$, the
function $y \rightarrow \tau_x f(y)$ belongs to
$\mathcal{D}(\mathbb{R}^d)$  and we have $$\forall\; y \in
\mathbb{R}^d,\; \mathcal{F}_D(\tau_x f)(y) = K(ix, y)
\mathcal{F}_D(f)(y).\eqno{(5.11)} $$ \hspace*{5mm} $\nu$i) For $f$
in $\mathcal{D}(\mathbb{R}^d)$ we have $$\forall\; x \in
\mathbb{R}^d,\; \int_{ \mathbb{R}^d}\tau_x f(y) \omega_k(y)dy =
\int_{ \mathbb{R}^d} f(y) \omega_k(y)dy.\eqno{(5.12)}$$
\hspace*{5mm} At the moment an explicit formula for the Dunkl
 translation operators is known only in the following two cases.(See [20,22]). \\
\underline{1$^{st}$ cas }: $d = 1$ and $W = \mathbb{Z}_2$. \\ For
all $f$ in $C(\mathbb{R})$ we have $$
 \forall \, x \in \mathbb{R}, \tau_{y}f(x)  =
 \frac{1}{2}\int_{-1}^{1}f(\sqrt{x^2 + y^2 -2xyt})
  (1+\frac{x-y}{\sqrt{x^2 + y^2 -2xyt}})\Phi_k(t)dt$$ $$
   \hspace*{3cm} +   \frac{1}{2}\int_{-1}^{1}f(-\sqrt{x^2 + y^2 -2xyt})
  (1-\frac{x-y}{\sqrt{x^2 + y^2 -2xyt}})\Phi_k(t)dt,
  \eqno{(5.13)}
$$ where $$\Phi_k(t) = \frac{\Gamma(k+\frac{1}{2})}{\sqrt{\pi}%
\Gamma(k)} (1+t)(1-t^2)^{k-1}.\eqno{(5.14)}$$
 \underline{2$^{nd}$
cas }: For all $f$ in $C( \mathbb{R}^d)$ radial we have $$ \forall
\, x \in \mathbb{R}^d, \; \tau_{y}f(x) = V_k [f_0 (\sqrt{||x||^2 +
||y||^2 +2 \langle x,.\rangle })](y),\eqno{(5.15)}$$ with $f_0$
the function on $[0,+\infty[$ given by $$f(x) =
f_0(||x||).\eqno{(5.16)}.$$ \hspace*{5mm}The Dunkl convolution
product of  $f$ and $g$ in $\mathcal{D}(\mathbb{R}^d) $ is the
functions $f*_{D}g$ defined by $$ \forall \, x \in \mathbb{R}^d,
\;
f*_{D}g(x)=\int_{\mathbb{R}^d}\tau_{x}f(-y)g(y)d\omega_{k}(y).\eqno{(5.17)}.$$
This convolution  is commutative and associative and admits  the
following properties\\ i) For $f$, $g$ in
$\mathcal{D}(\mathbb{R}^d)$( resp. $\mathcal{S}( \mathbb{R}^d)$)
the function $f*_{D}g$ belongs to $\mathcal{D}(\mathbb{R}^d)$(
resp. $\mathcal{S}( \mathbb{R}^d)$) and we have $$ \forall y \in
\mathbb{R}^d, \; {\cal F}_D (f*_{D}g)(y) = {\cal F}_D (f)(y){\cal
F}_D (g)(y).\eqno{(5.18)}.$$ \hspace*{5mm}ii)  For $f$, $g$ in
$\mathcal{D}(\mathbb{R}^d)$( resp. $\mathcal{S}( \mathbb{R}^d)$)
we have $$ \,^{t}V_k( f*_{D}g) = \,^{t}V_k (f)* \,^{t}V_k
(g)\eqno{(5.19)}$$ where $*$ is the classical convolution product
of functions on $ \mathbb{R}^d$.
\subsection{Dunkl Convolution product of distributions }
{\bf D\'efinition  5.1.} The Dunkl Convolution product of a
distribution $S$  in $\mathcal{D}'(\mathbb{R}^d)$ and a function
$\varphi$ in $\mathcal{D}(\mathbb{R}^d)$  is the function $S\ast_D
\varphi$ defined by $$\forall\; x \in \mathbb{R}^d,\;
S\ast_D\varphi(x) = \langle S_y, \tau_{-y}
\varphi(x)\rangle.\eqno{(5.20)}$$ {\bf Remark 5.1}

If $S = T_{f\omega_k}$ is the distribution  in
$\mathcal{D}'(\mathbb{R}^d)$ given by the  function $f\omega_k$
with $f$ in $ C(\mathbb{R}^d)$, we have  $$S\ast_D \varphi = f
\ast_D \varphi .\eqno{(5.21)}$$ {\bf Theorem 5.1.} The function $S
\ast_D \varphi$ is of class $C^{\infty}$ on $\mathbb{R}^d$ and we
have $$T^\mu (S\ast_D \varphi) = S \ast_D(T^\mu \varphi) = (T^\mu
S) \ast_D \varphi,\eqno{(5.22)}$$ where $$T^\mu =
T_1^{\mu_1}o\;T_2^{\mu_2}o\cdots\;o T_d^{\mu_d}, \mbox{ with } \mu
= (\mu_1, \mu_2,\cdots,\mu_d) \in \mathbb{N}^d,$$  and $T_j$, $j =
1,2,...,d$, the Dunkl operator defined on
$\mathcal{D}'(\mathbb{R}^d)$ by $$ \langle T_j S,\varphi \rangle =
- \langle S,T_j \varphi\rangle, \; \varphi \in \mathcal{D}(
\mathbb{R}^d). \eqno{(5.23)}$$ {\bf Proof}\begin{itemize}
\item[i)]We shall prove first that the function $S\ast_D \varphi$
is continuous on $ \mathbb{R}^d$. Let $x^0 \in \mathbb{R}^d$ and
$B(x^0 , r)$ the closed ball of center $x^0$ and radius $r > 0$.
We consider $\varphi$ in $\mathcal{D}(\mathbb{R}^d)$ such that
$\mbox{supp} \varphi \subset B(o,a)$, $a > 0$. From the relation
(5.4) and the fact that the support of $\mu_x$ is contained in the
closed ball of center $o$ and radius $||x||$, the relation (3.10)
implies that $$\forall x \in B(x^0,r), \; \mbox{supp}
\tau_{-y}\varphi(x) \subset B(o,a+r+||x^0||).$$ We put $K =
B(o,a+r+||x^0||)$ and $\Phi(x,y) = \tau_{-y}\varphi(x).$ \\ As the
distribution $S$ belongs to $\mathcal{D'}(\mathbb{R}^d)$ then
there exist a seminorm $p_n$ and a positive constant $C$ such that
for all $\theta$ in $\mathcal{D}(\mathbb{R}^d)$ with $\mbox{supp}
\theta \subset K,$ we have $$|\langle S,\theta\rangle| < C p_n
(\theta).\eqno{(5.24)}$$On the other hand as the function $\Phi$
is of class $C^{\infty}$ on $\mathbb{R}^d \times \mathbb{R}^d$,
then from (5.3) for all $\alpha \in \mathbb{N}^d$ the function
$D^\alpha \Phi(x,y)$ is continuous on $\mathbb{R}^d\times K$.
Thus: $\forall \varepsilon > 0, \, \exists \, B(x^0,r_0)$, $r_0 <
r$, \, $\forall \, x \in B(x^0,r_0) $ $\Longrightarrow$
$p_n(\Phi(x,.)-\Phi(x^0,.)) < \frac{\varepsilon}{C}$.\\ This
relation and (5.24) imply $$\forall \, x \in B(x^0,r_0), \; |S*_D
\varphi(x) - S*_D \varphi(x^0)| \leq C p_n(\Phi(x,.)-\Phi(x^0,.))
< \varepsilon. $$ The function $S*_D \varphi$ is then continuous
at $x^0$, and thus it is continuous on $\mathbb{R}^d$.
\item[ii)]We shall prove now that the function $S*_D
\varphi$ admits a partial derivative $\frac{\partial}{\partial
x_j}S*_D \varphi$ at $x^0 \in \mathbb{R}^d$ and we have
$$\frac{\partial}{\partial x_j}S*_D \varphi(x^0) = \langle
S_y,\frac{\partial}{\partial x_j}\tau_{-y}
\varphi(x^0)\rangle.\eqno{(5.25)}$$ \hspace*{5mm} Let $h \in
\mathbb{R}^d\backslash\{0\}$, by applying the Taylor formula we
obtain $$\begin{array}{lll}
  \Phi(x_1^0,...,x_j^0 + h,...,x_d^0,y)  & = & \Phi(x^0,y) +
  h\frac{\partial}{\partial x_j}\Phi(x^0,y) \\
   & + & h^2 \int_0^1 (1-t)\frac{\partial^2}{\partial x_j^2} \Phi(x_1^0,...,x_j^0 +
  th,...,x_d^0,y)dt.
\end{array}$$
Thus $$\frac{S \ast_D\varphi(x) - S\ast_D\varphi(x^0)}{h} -
\langle S_y,\frac{\partial}{\partial x_j}\Phi(x^0,y)\rangle = h
\langle S_y,R(x^0,y,h)\rangle, $$ where $$R(x^0,y,h) = \int_0^1
(1-t)\frac{\partial^2}{\partial x_j^2} \Phi(x_1^0,...,x_j^0 +
   th,...,x_d^0,y)dt.\eqno{(5.26)}$$
   It suffices to show that $\langle S_y,R(x^0,y,h)\rangle$
   remains bounded when $h$ tends to zero. We put $$M = \sup_{\begin{array}{lll}
  \hspace*{9mm} |\alpha|\leq n\\ (x,y) \in B(x^0,r)\times K\end{array}}
   |D_y^\alpha \frac{\partial^2}{\partial x_j^2}\Phi(x,y)|.$$ From
   (5.26) we deduce that $$\forall y \in K, \; |D_y^\alpha R(x^0,y,h)| \leq
   M.$$ From (5.24) we deduce that  $$|\langle S_y,R(x^0,y,h)\rangle| \leq
   CM.$$ Thus the  function $S \ast_D\varphi$ admits the partial
   derivative $\frac{\partial}{\partial x_j}S \ast_D\varphi(x^0)$
   at $x^0 \in \mathbb{R}^d$ and we have (5.25).\\ These result
   are true on $\mathbb{R}^d$ and in particular we have
   $$\forall \, x \in \mathbb{R}^d, \; \frac{\partial}{\partial x_j}S
   \ast_D\varphi(x)= \langle S_y,\frac{\partial}{\partial x_j} \tau_{-y}
   \varphi(x)\rangle.\eqno{(5.27)}
   $$
By applying the i) to this partial derivative we deduce that it is
continuous on $\mathbb{R}^d$.\\ Similar proofs as for i) and ii)
show that the function $S \ast_D\varphi$ admits continuous partial
derivatives of all order with respect to all variables. Then the
function $S \ast_D\varphi$ is of class $C^\infty$ on
$\mathbb{R}^d$. On the other hand using the definition of Dunkl
operator $T_j$ and the relations (5.27),(5.9) we obtain
$$\hspace*{-4cm}\forall \, x \in \mathbb{R}^d, \; T_j (S\ast_D
\varphi)(x) = \langle S_y,T_j \tau_{-y}
   \varphi(x)\rangle$$
$$\hspace*{3,8cm} = \langle S_y, \tau_{-y}(T_j
   \varphi)(x)\rangle = S\ast_D (T_j \varphi)(x).\eqno{(5.28)}
   $$
By iteration we get $$\forall \, x \in \mathbb{R}^d, \; T^\mu
(S\ast_D \varphi)(x) = S\ast_D (T^\mu \varphi)(x).$$ On the other
hand from (5.28) and (5.10) we have $$\begin{array}{lll} \forall\,
x \in \mathbb{R}^d, \; T_j (S\ast_D \varphi)(x)  & = &  \langle
S_y,(T_j)_y ( \tau_{-y}
   \varphi(x))\rangle\\
   & = &  \langle T_j S_y, \tau_{-y}
   \varphi(x)\rangle\\
   & = &  (T_j
S)\ast_D \varphi(x).
\end{array}$$
By applying this relation to the other Dunkl operators and their
composite we obtain $$\forall \, x \in \mathbb{R}^d, \;
T^\mu(S\ast_D \varphi)(x) = (T^\mu S)\ast_D \varphi(x).$$ This
completes the proof of the theorem.
\end{itemize} \noindent{\bf{Remark 5.2.}}\\ \hspace*{5mm} We have
$$\forall x \in \mathbb{R}^d, \; \delta\ast_D \varphi(x) = \langle
\delta, \tau_{-y}
   \varphi(x)\rangle = \varphi(x).$$
\subsection{Tensoriel product of distributions (see [21][3])}
\noindent{\bf{Theorem 5.2.}} Let $S, \mathcal{U} $ be two
distributions in $ \mathcal{D'}(\mathbb{R}^d).$ Then  \\
\hspace*{5mm} i) There exists a unique distribution in
$\mathcal{D'}(\mathbb{R}^d\times \mathbb{R}^d)$ such that for all
$\varphi,\psi$ in $\mathcal{D}(\mathbb{R}^d)$ we have $$\langle
W,\varphi\otimes\psi \rangle = \langle S,\varphi\rangle \langle
\mathcal{U},\psi\rangle.\eqno{(5.29)}$$ The distribution $W$ is
the tensoriel product of the distributions $S$ and $\mathcal{U}$
and it is denoted by $S\otimes \mathcal{U}$. \\\hspace*{5mm} ii)
For all $\phi$ in $\mathcal{D}(\mathbb{R}^d\times\mathbb{R}^d)$ we
have $$\langle S\otimes \mathcal{U},\phi(x,y) \rangle = \langle
S_x,\langle \mathcal{U}_y, \phi(x,y)\rangle\rangle = \langle
\mathcal{U}_y,\langle S_x, \phi(x,y)\rangle\rangle.\eqno{(5.30)}$$
The tensoriel product of distributions satisfies the following
properties.\\ \hspace*{5mm} i) Let $S,\mathcal{U},\mathcal{V}$ be
in $ \mathcal{D'}(\mathbb{R}^d)$. There exists a unique
distribution $S\otimes \mathcal{U}\otimes \mathcal{V}$ in
$\mathcal{D'}(\mathbb{R}^d\times\mathbb{R}^d\times)$ such that for
all $\varphi,\psi,\theta$ in $\mathcal{D}(\mathbb{R}^d)$ we have
$$\langle S\otimes \mathcal{U}\otimes
\mathcal{V},\varphi\otimes\psi\otimes\theta \rangle = \langle
S,\varphi\rangle \langle \mathcal{U},\psi\rangle\langle
\mathcal{V},\theta\rangle.\eqno{(5.31)}$$ We deduce that
$$(S\otimes \mathcal{U})\otimes \mathcal{V} = S\otimes
(\mathcal{U}\otimes \mathcal{V}).\eqno{(5.32)}$$
 \hspace*{5mm} ii) For all  $S, \mathcal{U}$ in $
\mathcal{D'}(\mathbb{R}^d)$ we have $$D^\alpha(S\otimes
\mathcal{U}) = (D^\alpha S)\otimes \mathcal{U}.\eqno{(5.33)}$$
 \hspace*{5mm} iii) For all  $S,\mathcal{U}$ in $
\mathcal{D'}(\mathbb{R}^d)$ we have $$\mbox{supp}(S\otimes
\mathcal{U}) = (\mbox{supp} S)\times (\mbox{supp}
\mathcal{U}).\eqno{(5.34)}$$
\subsection{Dunkl convolution product of distributions}
\hspace*{5mm} To define the Dunkl convolution product of the
distributions $S$ and $\mathcal{U}$ in  $
\mathcal{D'}(\mathbb{R}^d)$ we must consider the expression
$$\langle S_x\otimes \mathcal{U}_y, \tau_x \varphi(y)\rangle, \;
\varphi \in \mathcal{D}(\mathbb{R}^d).\eqno{(5.35)}$$ But the
function $(x,y) \to \tau_x \varphi(y)$ defined on
$\mathcal{D}(\mathbb{R}^{d}\times \mathbb{R}^d)$ is not with
compact support. Then the expression (5.35) is not well defined in
the general case. It will have a sense   if the set $$\mbox{supp}
(\tau_x \varphi(y)) \cap \mbox{supp}(S\otimes \mathcal{U})$$ is
compact for all $\varphi$ in $\mathcal{D}(\mathbb{R}^d)$. We say
in this case that the supports of $S$ and $\mathcal{U}$ satisfy
the " supports condition ".\\ \noindent{\bf{Definition 5.2.}} Let
$S$ and $\mathcal{U}$ be two distributions in
$\mathcal{D'}(\mathbb{R}^d)$ such that their supports satisfy the
" supports condition ". The Dunkl convolution product of $S$ and
$\mathcal{U}$ is the distribution $S*_D \mathcal{U}$ in
$\mathcal{D'}(\mathbb{R}^d)$ defined by $$\langle S*_D
\mathcal{U}, \varphi\rangle = \langle S_x\otimes \mathcal{U}_y,
\tau_x \varphi(y)\rangle.\eqno{(5.36)}$$ \noindent{\bf{Remark
5.3.}}
\\ \hspace*{5mm} Using (5.30) the relation (5.36) can also be
written in the form $$\langle S*_D \mathcal{U}, \varphi\rangle =
\langle S_x,\langle \mathcal{U}_y, \tau_x \varphi(y)\rangle\rangle
= \langle \mathcal{U}_y,\langle S_x, \tau_x
\varphi(y)\rangle\rangle .\eqno{(5.37)}$$
 \noindent{\bf{Proposition 5.1.}} Let $S$ be in
 $\mathcal{E'}(\mathbb{R}^d)$ and $\mathcal{U} $ in
 $\mathcal{D'}(\mathbb{R}^d)$. Then the supports of these
 distributions  satisfies the " supports condition.
"\\

\hspace*{5mm} The Dunkl convolution product of distributions is
commutative and associative and  satisfies the following
properties.\\ \hspace*{5mm}i) Let $S$  and $\mathcal{U} $ be in
 $\mathcal{D'}(\mathbb{R}^d)$ such that their support satisfy the " supports condition.
" Then for all $\mu \in \mathbb{N}^d$ we have $$T^\mu(S*_D
\mathcal{U}) = (T^\mu S)*_D \mathcal{U} = S*_D (T^\mu
\mathcal{U}).\eqno{(5.38)}$$\hspace*{5mm}ii) Let $S$ be  in
 $\mathcal{D'}(\mathbb{R}^d)$ and $f \in  \mathcal{D}(\mathbb{R}^d)$.
 Then  we have $$S*_D T_{f\omega_k} =T_{(S*_D f)\omega_k}.\eqno{(5.39)}$$
where $T_{f\omega_k}$ is the distribution in
$\mathcal{D'}(\mathbb{R}^d)$ given by the function
$f\omega_k$.\\\noindent{\bf{Remark 5.4}}\\ \hspace*{5mm} Let $S$
and $U$ be two distributions in $ \mathcal{E'}( \mathbb{R}^d)$,
with $\mbox{supp}\, S \subset B(o,a)$, $a > 0$, and $\mbox{supp}\,
U \subset B(o,b)$, $b > 0$. Then the distribution $S*_D U$ belongs
to  $ \mathcal{E'}( \mathbb{R}^d)$ and we have $\mbox{supp}\, S*_D
U \subset B(o,a+b).$\\ {\bf Proposition 5.2.} Let $S$  and $U$ be
two distributions in $\mathcal{E}'(\mathbb{R}^d)$. Then  we have
$$\,^{t}V_k ( S \ast_D U) = \,^{t}V_k(S)* \,\,^{t}V_k(U).$$ where
$*$ is the classical convolution product of distributions on
$\mathbb{R}^d$.\\{\bf Proof}\\ \hspace*{5mm} From (3.14) for
$\varphi$ in $\mathcal{E}(\mathbb{R}^d)$ we have
$$\langle\,\,^{t}V_k ( S \ast_D \mathcal{U}),\varphi\rangle =
\langle  S \ast_D U,V_k(\varphi)\rangle.$$ By using (5.37) and
(5.1) we obtain $$\begin{array}{lll}
   \langle\,\,^{t}V_k ( S \ast_D
\mathcal{U}),\varphi\rangle & = & \langle S_x,\langle U_y,
\tau_x(V_k(\varphi))(y)\rangle\rangle
\\
   & = & \langle S_x,\langle U_y,
(V_k)_x((V_k)_y[\varphi(x+y)]\rangle\rangle.
\end{array}$$By applying (3.14) we get $$\begin{array}{lll}
   \langle\,\,^{t}V_k ( S \ast_D
\mathcal{U}),\varphi\rangle & = & \langle \,\,^{t}V_k
(S)_x,\langle \,\,^{t}V_k (U)_y, \varphi(x+y)\rangle\rangle\\ &=&
\langle \,^{t}V_k (S)*\; \,^{t}V_k (U),\varphi\rangle.
\end{array}$$
Thus $$\,^{t}V_k ( S \ast_D U) = \,^{t}V_k(S)* \; \,^{t}V_k(U).$$
\subsection{Dunkl Convolution product of  tempered distributions}
\hspace*{5mm} The results of this subsection are proved in [26].\\
{\bf D\'efinition 5.3} Let $S$ be in $\mathcal{S}'(\mathbb{R}^d)$
and $\varphi$ in $\mathcal{S}(\mathbb{R}^d)$. The Dunkl
Convolution product of $S$ and $\varphi$    is the function
$S\ast_D \varphi$ defined by $$\forall\; x \in \mathbb{R}^d,\;
S\ast_D\varphi(x) = \langle S_y, \tau_{-y}
\varphi(x)\rangle.\eqno{(5.40)}$$
 {\bf Proposition 5.3.} For $S$
in $\mathcal{S}'(\mathbb{R}^d)$ and  $\varphi$ in
$\mathcal{S}(\mathbb{R}^d)$ the function $S\ast_D \varphi$ belongs
to $\mathcal{E}(\mathbb{R}^d)$ and we have $$T^\mu (S\ast_D
\varphi) = S \ast_D(T^\mu \varphi) = (T^\mu S) \ast_D
\varphi,\eqno{(5.41)}$$ where $$T^\mu =
T_1^{\mu_1}o\;T_2^{\mu_2}o\cdots\;o T_d^{\mu_d}, \mbox{ with } \mu
= (\mu_1, \mu_2,\cdots,\mu_d) \in \mathbb{N}^d.$$
\subsection{Dunkl transform of distributions}
{\bf Definition 5.4} \begin{itemize}
\item[i)] The Dunkl transform of a distribution $S$ in
$\mathcal{S}'(\mathbb{R}^d)$ is defined by $$\langle
\mathcal{F}_D(S), \psi\rangle = \langle S,
\mathcal{F}_D(\psi)\rangle, \psi \in
\mathcal{S}(\mathbb{R}^d).\eqno{(5.42)}$$
\item[ii)] We define the Dunkl transform of a distribution $S$ in
$\mathcal{E}'(\mathbb{R}^d)$ by $$\forall\; y \in \mathbb{R}^d,\;
\mathcal{F}_D(S)(y) = \langle S_x, K(-iy,x)\rangle.\eqno{(5.43)}
$$
\end{itemize}
{\bf Remarks 5.5}\begin{itemize}
\item[i)]

When the distribution $S$ in $\mathcal{E}'(\mathbb{R}^d)$ is given
by the function $g \omega_k$ with $g$ in
$\mathcal{D}(\mathbb{R}^d)$, and denoted by $T_{g \omega_k}$, the
relation (5.43) coincides with (4.1).
\item[ii)] From (3.14) and (2.16) the relation (5.43) can also be
written in the form $$\forall \, y \in \mathbb{R}^d, \; {\cal F}_D
(S)(y) = {\cal F}\circ\, \,^{t}V_k(S)(y), \eqno{(5.44)}$$ where
${\cal F}$ is the classical Fourier transform of distributions in
$ \mathcal{E'}( \mathbb{R}^d)$ given by $$\forall \, y \in
\mathbb{R}^d, \; {\cal F}(U)(y) = \langle U,e^{-i\langle
.,y\rangle}\rangle. \eqno{(5.45)}$$
\end{itemize}
 {\bf Notation.}\\ \hspace*{5mm}
We denote by $\mathcal{H}(\mathbb{C}^d)$ the space of entire
functions on $\mathbb{C}^d$ which are rapidly increasing and of
exponential type. We have $$\mathcal{H}(\mathbb{C}^d) = \bigcup_{
a \geq 0} \mathcal{H}_a(\mathbb{C}^d) $$ where
$\mathcal{H}_a(\mathbb{C}^d)$ is the space of entire functions
$\Psi$ on $\mathbb{C}^d$ satisfying $$\exists \, N \in \mathbb{N},
\; \sup_{z \in \mathbb{C}^d}(1+||z||^2)^{-N}|\Psi(z)|e^{-a||Im
z||} < +\infty. $$We topology this space with the classical
topology.\vspace{5mm}

The following theorem is given in [25] p.27.\\ {\bf Theorem 5.3.}
The transform $\mathcal{F}_D$ is a topological isomorphism from
\begin{itemize}
\item[i)] $\mathcal{S}'(\mathbb{R}^d)$ onto itself.
\item[ii)] $\mathcal{E}'(\mathbb{R}^d)$ onto
$\mathcal{H}(\mathbb{C}^d)$.
\end{itemize}
{\bf Theorem 5.4.} Let $S$ be in $\mathcal{S}'(\mathbb{R}^d)$ and
$\varphi$ in $\mathcal{S}(\mathbb{R}^d)$. Then the distribution on
$\mathbb{R}^d$ given by $(S\ast_D \varphi)\omega_k$ belongs to
$\mathcal{S}'(\mathbb{R}^d)$ and we have
$$\mathcal{F}_D(T_{(S\ast_D \varphi)\omega_k})=
\mathcal{F}_D(\varphi)\mathcal{F}_D(S).\eqno{(5.46)}$$
\hspace*{5mm}We consider the radial positive function $\varphi$ in
$\mathcal{D}(\mathbb{R}^d)$, with \mbox{support} in the closed
ball of center $0$ and radius 1, satisfying
$$\int_{\mathbb{R}^d}\varphi(x) \omega_k(x)dx = 1,$$ and $\phi$
the function on $[0, + \infty[$ given by $$\varphi(x) =
\phi(\|x\|) = \phi(r), \mbox{ with } r = \|x\|. $$ For
$\varepsilon \in ]0,1]$, we denote by $\varphi_\varepsilon$ the
function on $\mathbb{R}^d$ defined by $$\forall\; x \in
\mathbb{R}^d,\;\; \varphi_\varepsilon(x) =
\frac{1}{\varepsilon^{2\gamma+d}}
\phi(\frac{\|x\|}{\varepsilon}).\eqno{(5.47)}$$ {\bf Theorem 5.5.}
Let $S$ be in $\mathcal{S}'(\mathbb{R}^d)$. We have
$$\lim_{\varepsilon \rightarrow 0}S \ast_D \varphi_\varepsilon =
S,\eqno{(5.48)}$$ the limit is in $\mathcal{S}'(\mathbb{R}^d)$.\\
 {\bf Definition
5.5. } Let $S$ be in $\mathcal{S}'(\mathbb{R}^d)$ and
$\mathcal{V}$ in $\mathcal{E}'(\mathbb{R}^d)$. The Dunkl
convolution product of $S$ and $\mathcal{V}$ is the distribution
$S \ast_D \mathcal{V}$ on $\mathbb{R}^d$ defined by $$\langle S
\ast_D \mathcal{V},\psi\rangle = \langle S_{x}, \langle
\mathcal{V}_{y}, \tau_x\psi(y) \rangle\rangle,\; \psi \in
\mathcal{D}(\mathbb{R}^d).\eqno{(5.49)} $$ {\bf Remark 5.6.}

The relation (5.49) can also be written in the form $$\langle S
\ast_D \mathcal{V},\psi\rangle = \langle S, \check{\mathcal{V}}
\ast_D \psi\rangle, \; \psi \in \mathcal{D}(
\mathbb{R}^d),\eqno{(5.50)} $$ with $\check{\mathcal{V}}$ the
distribution in $\mathcal{E}'(\mathbb{R}^d)$ defined by $$\langle
\check{\mathcal{V}}, f\rangle = \langle
\mathcal{V},\check{f}\rangle,\; f\in \mathcal{E}(\mathbb{R}^d), $$
where $ \check{f}$ given by $$\forall \, x \in \mathbb{R}^d, \,
\check{f}(x) = f(-x).$$
 {\bf Theorem 5.6.} Let $S$ be in
$\mathcal{S}'(\mathbb{R}^d)$ and $\mathcal{V}$ in
$\mathcal{E}'(\mathbb{R}^d)$. Then the distribution $S \ast_D
\mathcal{V}$ belongs to $\mathcal{S}'(\mathbb{R}^d)$ and we have
$$\mathcal{F}_D(S \ast_D \mathcal{V}) =
\mathcal{F}_D(\mathcal{V}).\mathcal{F}_D(S).\eqno{(5.51)}$$ {\bf
Proof}

We deduce the result from (5.50),(5.39) and Theorem 5.4.\\ {\bf
Definition 5.6.} We define the dual Dunkl intertwining operator
$\,^{t}V_k $ on $ \mathcal{S'}( \mathbb{R}^d)$ by $$ \,^{t}V_k (S)
= {\cal F}^{-1}\circ{\cal F}_D(S). \eqno{(5.52)}$$
 {\bf Theorem
5.7.} \begin{itemize} \item [i)] The operator $ \,^{t}V_k $ is a
topological isomorphism from $ \mathcal{S'}( \mathbb{R}^d)$ onto
itself.
\item [ii)]  Let $S$ be in $ \mathcal{S'}( \mathbb{R}^d)$
and $\mathcal{U}$ in $ \mathcal{E'}( \mathbb{R}^d)$. Then we have
$$\,^{t}V_k ( S \ast_D \mathcal{U}) = \,^{t}V_k(S)*
\,\,^{t}V_k(\mathcal{U}).\eqno{(5.53)}$$ where $*$ is the
classical convolution product of tempered  distributions on $
\mathbb{R}^d$.
\end{itemize}{\bf Proof} \\ \hspace*{5mm}i) We deduce the result
from Theorem 5.3 and the properties of the classical Fourier
transform of tempered distributions on $ \mathbb{R}^d$.\\
\hspace*{5mm} ii) From Theorem 5.6 we have $$\mathcal{F}_D(S
\ast_D \mathcal{U}) =
\mathcal{F}_D(S).\mathcal{F}_D(\mathcal{U}).$$ Using (5.52) we
obtain $$\mathcal{F}\circ\,\,^{t}V_k(S \ast_D \mathcal{U}) =
\mathcal{F}\circ\,\,^{t}V_k(S).\mathcal{F}\circ\,\,^{t}V_k(\mathcal{U}).$$
Thus $$\mathcal{F}(\,^{t}V_k(S \ast_D \mathcal{U})) =
\mathcal{F}(\,^{t}V_k(S)*\,\,^{t}V_k(\mathcal{U})).$$ We deduce
(5.53) from this relation and the injectivity of the transform
${\cal F}$ on $ \mathcal{S'}( \mathbb{R}^d)$ .\\{\bf Remark 5.7.}
\\ \hspace*{5mm} When the distribution $S$ is in $ \mathcal{E'}(
\mathbb{R}^d)$ another proof of the relation (5.54) has been given
in Proposition 5.2.
\section{Hypoelliptic Dunkl convolution equations in the   space of distributions}

Let $S$ be in ${\cal E'}( \mathbb{R}^d)$. In this section we study
convolution equations of the form $$S *_D \mathcal{U} =
\mathcal{V}, \eqno{(6.1)}$$ where $\mathcal{U}$ and $\mathcal{V}$
are distributions in $ \mathcal{D'}( \mathbb{R}^d)$.\\
\hspace*{5mm} We say that the equation (6.1) is hypoelliptic if
all solution $\mathcal{U}$ is given by a function $f\omega_k$ with
$f$ in $\mathcal{E}( \mathbb{R}^d)$ whenever $\mathcal{V}$ is
given by a function $g\omega_k$ with $g$ in $\mathcal{E}(
\mathbb{R}^d)$.\\ \hspace*{5mm} When (6.1) is hypoelliptic we say
also that the distribution $S$ is hypoelliptic.\\ \hspace*{5mm}
The main result of this section is the characterization of
hypoelliptic Dunkl convolution equations in terms of their Dunkl
transform.\\ \hspace*{5mm}We say that the distribution $S$ in
${\cal E'}(\mathbb{R}^d)$ satisfies the $H$-property if \\
\hspace*{5mm} i) There exists $A, M > 0$ such that $|{\cal
F}_D(S)(x)| \geq ||x||^{-A}$ for all $||x|| \geq
M$.\\\hspace*{8mm}ii) $\lim_{||z|| \to \infty, \, z \in {\cal Z}}
\frac{||Im z||}{Log ||z||} = \infty,$ where ${\cal Z} = \{z \in
\mathbb{C}^d, \; {\cal F}_D(S)(z) = 0.\}$, with $||z||^2 =
\sum_{j=1}^d (Re z_j)^2 + (Im z_j)^2.$\\ \noindent{
\bf{Proposition 6.1.}} Let $S$ be in $ {\cal E'}(\mathbb{R}^d)$.
If $S$ is hypoelliptic  then $S$ satisfies the i) of the
$H$-property.\\ \hspace*{5mm} To prove this proposition we need
the following Lemma.
\\{\bf{Lemma 6.1.}} Let $\phi$ be a positive function in $D(
\mathbb{R}^d)$ such that $\phi(0) = 1$ and which is even for $d =
1$ and radial for $d \geq 2$. Then there exist positive constants
$C$ and $X$ such that for $||x|| \geq X$ we have $${\cal
F}_D(K(ix,.)\phi)(x) \geq
\frac{C}{||x||^{2\gamma+d}}.\eqno{(6.2)}$$ \noindent{\bf{Proof}}
\\ \hspace*{5mm} i) We suppose that $d \geq 2$. \\We have $$\forall x \in \mathbb{R}^d, \; {\cal
F}_D(K(ix,.)\phi)(x) = \int_{ \mathbb{R}^d}|K(ix,t)|^2
\phi(t)\omega_k(t)dt.$$ As $\phi$ is radial then there exists a
function $\varphi$ on $[0,+\infty[$ such that $$\phi(t) =
\varphi(||t||) = \varphi(r), \; with \, r = ||t||.$$ By using
polar coordinates we obtain $$\forall x \in \mathbb{R}^d, \; {\cal
F}_D(K(ix,.)\phi)(x) = \int_0^{\infty}(\int_{
S^{d-1}}|K(ix,r\sigma)|^2
\omega_k(\sigma)d\sigma)\varphi(r)r^{2\gamma+d-1}dr.$$ As the
function $\varphi$ is positive, then for all $x \in
\mathbb{R}^d\backslash\{0\}$ we have $$\begin{array}{lll} {\cal
F}_D(K(ix,.)\phi)(x) &\geq &\int_0^{\frac{1}{||x||}}(\int_{
S^{d-1}}|K(i||x||\beta,r\sigma)|^2
\omega_k(\sigma)d\sigma)\varphi(r)r^{2\gamma+d-1}dr\\ &\geq&
\int_0^{\frac{1}{||x||}}(\int_{ S^{d-1}}|K(i\beta,r||x||\sigma)|^2
\omega_k(\sigma)d\sigma)\varphi(r)r^{2\gamma+d-1}dr.
\end{array}$$
By the change of variables $u = r||x||$, we obtain $${\cal
F}_D(K(ix,.)\phi)(x) \geq
\frac{1}{||x||^{2\gamma+d}}\int_0^{1}(\int_{
S^{d-1}}|K(i\beta,u\sigma)|^2
\omega_k(\sigma)d\sigma)\varphi(\frac{u}{||x||})u^{2\gamma+d-1}du.\eqno{(6.3)}$$
We denote by $I_\beta(||x||)$ the integrals of the second member.
From the properties of the function $\varphi$ we deduce that there
exists $X > 0$ such that  for all $u \in [0,1]$ and $||x|| \geq
X$, we have $\varphi(\frac{u}{||u||}) \geq \frac{1}{2}$. Then
$$I_\beta(||x||) \geq \frac{1}{2} \int_0^{1}(\int_{
S^{d-1}}|K(i\beta,u\sigma)|^2
\omega_k(\sigma)d\sigma)u^{2\gamma+d-1}du.$$ As the second member
is continuous on $S^{d-1}$ with respect to the variable $\beta$,
then for $||x|| \geq X$, we have $$I_\beta(||x||)  \geq
\frac{1}{2}\min_{\beta \in S^{d-1}} \int_0^{1}(\int_{
S^{d-1}}|K(i\beta,u\sigma)|^2
\omega_k(\sigma)d\sigma)u^{2\gamma+d-1}du,$$ and there exists
$\beta_0 \in S^{d-1}$ such that for $||x|| \geq X$, we have
$$I_\beta(||x||)  \geq \frac{1}{2} \int_0^{1}(\int_{
S^{d-1}}|K(i\beta_0,u\sigma)|^2
\omega_k(\sigma)d\sigma)u^{2\gamma+d-1}du.\eqno{(6.4)}$$ As the
function $$u \to \int_{ S^{d-1}}|K(i\beta_0,u\sigma)|^2
\omega_k(\sigma)d\sigma$$ is continuous on $[0,1]$, then if
$$\int_0^{1}(\int_{ S^{d-1}}|K(i\beta_0,u\sigma)|^2
\omega_k(\sigma)d\sigma)u^{2\gamma+d-1}du = 0,$$ we deduce that
$$\forall \, u \in [0,1], \; \int_{
S^{d-1}}|K(i\beta_0,u\sigma)|^2 \omega_k(\sigma)d\sigma = 0. $$ By
taking $u = 0$, we obtain $$\int_{ S^{d-1}}
\omega_k(\sigma)d\sigma = 0,$$ which contradicts $$\int_{ S^{d-1}}
\omega_k(\sigma)d\sigma = d_k = \frac{2}{c_k
\Gamma(\gamma+\frac{d}{2})} .$$ Then $$\int_0^{1}(\int_{
S^{d-1}}|K(i\beta_0,u\sigma)|^2
\omega_k(\sigma)d\sigma)u^{2\gamma+d-1}du \neq 0.\eqno{(6.5)}$$ We
denote the first member of this relation by $2C$. By using the
relations (6.3),(6.4) and (6.5), we deduce that for $||x|| \geq X
$, we have $${\cal F}_D(K(ix,.)\phi)(x) \geq
\frac{C}{||x||^{2\gamma+d}}.$$ \hspace*{5mm} ii) We suppose that
$d = 1$.\\ The same proof as for i) gives the relation (6.2).\\
\noindent{\bf{Proof of Proposition 6.1. }}\\ \hspace*{5mm} We
assume  that the i) of the $H$-property does not hold. Then we can
find a sequence $(x_n)_{n \in \mathbb{N}} \subset \mathbb{R}^d$
such that $||x_n|| \geq 2^n$ and $$ \forall \, n \in \mathbb{N},
\; |{\cal F}_{D}(S)(x_n)| < ||x_n||^{-n}.\eqno{(6.6)}$$  We
consider the sequence $\{U_p\}_{p\in \mathbb{N}}$ of distributions
in $ \mathcal{D'} ( \mathbb{R}^d)$ given by $$U_p =
\sum_{n=0}^{p}T_{K(-ix_n,.)\omega_k},$$ where $
T_{K(-ix_n,.)\omega_k}$  the distribution in $ \mathcal{D'}(
\mathbb{R}^d)$ given by the function $ K(-ix_n,.)\omega_k.$\\
 Let $\varphi$ be in $D( \mathbb{R}^d)$. For all $p,q
\in \mathbb{N}$ with $p > q$, we have $$\langle U_p,\varphi\rangle
- \langle U_q,\varphi\rangle = \sum_{n=q}^{p}\langle
T_{K(-ix_n,.)\omega_k},\varphi\rangle,$$ $$\hspace*{2,3cm}=
\sum_{n=q}^{p}{\cal
   F}_{D}(\varphi)(x_n).\eqno{(6.7)}$$
But from Theorem 4.1 the function ${\cal
   F}_{D}(\varphi)$ is rapidly decreasing. Then there exists a
   positive constant $C$ such that $$\forall \, y \in \mathbb{R}^d, \; |{\cal
   F}_{D}(\varphi)(y)| \leq \frac{C}{1+||y||}.$$
Thus $$\forall \, n \in \mathbb{N}, \; |{\cal
   F}_{D}(\varphi)(x_n)| \leq \frac{C}{||x_n||}\leq  \frac{C}{2^n}.\eqno{(6.8)}$$
By applying this relation to (6.7) we obtain $$|\langle
U_p,\varphi\rangle - \langle U_q,\varphi\rangle| \leq C
\sum_{n=q}^{p}\frac{1}{2^n}\to 0, \; as \, q \to +\infty.$$ Then
$$\langle U_p,\varphi\rangle \to L(\varphi), \; as \, p \to
+\infty.$$ We deduce that $L$ is a distribution $U$ in $
\mathcal{D'}( \mathbb{R}^d)$ and $U_p$ converges to $U$ in $
\mathcal{D'}( \mathbb{R}^d)$ as $p$ tends to infinity. Thus $$U =
\sum_{n=0}^{\infty}T_{K(-ix_n,.)\omega_k}, \eqno{(6.9)}$$ and for
all $\varphi$ in $\mathcal{D}( \mathbb{R}^d)$ we have $$\langle
U,\varphi\rangle = \sum_{n=0}^{\infty}{\cal
   F}_{D}(\varphi)(x_n).
\eqno{(6.10)}$$ \hspace*{5mm} We shall prove now that the
distribution $S*_D U$ of $ \mathcal{D'}( \mathbb{R}^d)$ is given
by a function $f\omega_k$ with $f$ in $ \mathcal{E}(
\mathbb{R}^d).$\\ \hspace*{5mm} From (5.37),(6.6) and (5.11), for
all $\varphi$ in $ \mathcal{D}( \mathbb{R}^d)$ we have $$ \langle
S*_D U,\varphi \rangle = \langle S_y,\langle U_t,\tau_y \varphi(t)
\rangle\rangle$$
 $$ \langle S*_D U,\varphi \rangle = \langle
S_y,\sum_{n=0}^\infty K(iy,x_n){\cal
   F}_{D}(\varphi)(x_n)\rangle. $$
By applying Theorem 4.1 and Definition 5.4 ii) we obtain $$\langle
S*_D U,\varphi \rangle = \sum_{n=0}^\infty {\cal
   F}_{D}(\varphi)(x_n)\langle
S_y,K(iy,x_n)\rangle$$ $$\hspace*{9mm}= \sum_{n=0}^\infty {\cal
   F}_{D}(\varphi)(x_n){\cal
   F}_{D}(S)(-x_n).
\eqno{(6.11)}$$ This relation can also be written in the form
$$\langle S*_D U,\varphi \rangle = \sum_{n=0}^\infty {\cal
   F}_{D}(S)(-x_n)\int_{ \mathbb{R}^d} K(-it,x_n)\varphi(t) \omega_k(t) dt.$$
By using (2.14) and the fact that $\varphi$ belongs to $
\mathcal{D}( \mathbb{R}^d)$ and ${\cal
   F}_{D}(S)$ satisfies the relation (6.6), we can interchange the
   series and the integral and we obtain $$ \langle S*_D U,\varphi \rangle =
   \int_{ \mathbb{R}^d}[\sum_{n=0}^{\infty}{\cal
   F}_{D}(S)(-x_n) K(-it,x_n)]\varphi(t) \omega_k(t) dt.$$
Thus the distribution  $S*_D U$ is given by the function
$f\omega_k$, with $$ f(t) = \sum_{n=0}^{\infty}{\cal
   F}_{D}(S)(-x_n) K(-it,x_n).$$
\hspace*{5mm} Let $t_0 \in \mathbb{R}^d$ and $B(t_0,r)$ the open
ball of center $t_0$ and radius $r > 0$. By using (6.6) and (2.12)
we deduce that for all $\nu \in \mathbb{N}^d$ there exists a
positive constant $C$ such that $$\sup_{t\in
B(t_0,r)}|D^{\nu}({\cal
   F}_{D}(S)(-x_n) K(-it,x_n))| \leq C ||x_n||^{-n+|\nu|}.$$
As the series $ \sum_{n=0}^\infty||x_n||^{-n+|\nu|}$ converges, we
deduce that the function $f$ admits continuous partial derivatives
of all order on $B(t_0,r)$ and then on $ \mathbb{R}^d$. Thus $f$
belongs to $ \mathcal{E}( \mathbb{R}^d)$.\\
 \hspace*{5mm} In the
following we want to show that the distribution $U$ does not given
by  a function $g\omega_k$ with $g $ in $ \mathcal{E}(
\mathbb{R}^d)$. If not we take a positive  function $\phi$ in $
\mathcal{D}( \mathbb{R}^d)$ such that $\phi(0) = 1$ which is even
for $d = 1$
  and radial for $d \geq 2$ and ${\cal
   F}_{D}^{-1}(\phi)$ is non negative. We consider $\triangle_k$
   the Dunkl Laplacian given by $$\triangle_k = \sum_{j=1}^d T_j^2,$$
and $p \in \mathbb{N}$ such that $p > \frac{1}{2}(\gamma +
\frac{d}{2}).$ Using (6.10) and (2.8) we obtain for all $y \in
\mathbb{R}^d$:

 $$\hspace*{-3,7cm}\langle K(iy,.)\triangle_k^{2p} U,\phi\rangle
 =
 \langle \triangle_k^{2p} U,K(iy,.)\phi\rangle\eqno{(6.12)}$$ $$
 \begin{array}{lll}
  \hspace*{1cm} &=& \langle \triangle_k^{2p}T_{g \omega_k},K(iy,.)\phi\rangle\\ &=&
   \langle T_{(\triangle_k^{2p}g)\omega_k},K(iy,.)\phi)\rangle\\ &=&
   \int_{\mathbb{R}^d}K(iy,t)\phi(t) \triangle_k^{2p}g(t)\omega_k(t)dt.  \end{array}$$
Bt taking $y = x_j$, we have $$\langle K(iy,.)\triangle_k^{2p}
U,\phi\rangle = \int_{\mathbb{R}^d}K(ix_j,t)\phi(t)
\triangle_k^{2p}g(t)\omega_k(t)dt. $$ As the function $\phi
\triangle_k^{2p}g$ belongs to $L^1_k( \mathbb{R}^d)$, then from
the properties of the Dunkl transform we have $$\lim_{j \to
+\infty}\langle K(iy,.)\triangle_k^{2p} U,\phi\rangle =
0.\eqno{(6.13)}$$ On the other hand from (6.12) we have $$
\hspace*{-3,5cm}\langle K(iy,.)\triangle_k^{2p} U,\phi\rangle =
\langle U,\triangle_k^{2p}(K(iy,.)\phi)\rangle
$$$$\hspace*{0,8cm}= \sum_{n=0}^{\infty}{\cal
   F}_{D}(\triangle_k^{2p}(K(iy,.)\phi))(x_n).$$
    Thus
   $$\langle K(iy,.)\triangle_k^{2p} U,\phi\rangle = \sum_{n=0}^{\infty}||x_n||^{4p}{\cal
   F}_{D}(K(iy,.)\phi)(x_n).\eqno{(6.14)}$$
On the other hand  for all $ y \in \mathbb{R}^d$ we have $$\forall
\, z \in \mathbb{R}^d, \; {\cal
   F}_{D}(K(iy,.)\phi)(z) = \int_{ \mathbb{R}^d} K(iy,t)K(-it,z)\phi(t)\omega_k(t)dt
   .\eqno{(6.15)}$$ Thus from Theorem 4.1 and the relation (5.11)
   we obtain
$$\forall \, z \in \mathbb{R}^d, \; {\cal
   F}_{D}(K(iy,.)\phi)(z) = \frac{2^{2\gamma+d}}{c_k^2}\tau_{y}({\cal
   F}_{D}^{-1}(\phi))(z).$$ As the function ${\cal
   F}_{D}^{-1}(\phi)$ is  positive and even for $d =1$ and radial for $d \geq 2$,
    then for $d \geq 2$ we have $$ \forall \, t \in \mathbb{R}^d, \; {\cal
   F}_{D}^{-1}(\phi)(t) = F(||t||),$$ where $F$ is a positive
   function on $[0,+\infty[$.\\ Using (5.15),(5.
   13) and example 3.1 we have
   $$\forall \, z \in \mathbb{R}^d, \; {\cal
   F}_{D}(K(iy,.)\phi)(z) = \frac{2^{2\gamma+d}}{c_k^2}
   V_k[F(\sqrt{||y||^2 + ||z||^2 - 2\langle y,.\rangle})](z),
   \, \mbox{if} \, d \geq 2.$$
   $$\forall \, z \in \mathbb{R}, \; {\cal
   F}_{D}(K(iy,.)\phi)(z) = C \,
   V_k[({\cal
   F}_{D}^{-1}(\phi))(\sqrt{y^2 + z^2 - 2yz})](z)\, \mbox{if} \, d = 1,$$
   with $$C = 2^{2k+1}(\Gamma(k+\frac{1}{2}))^2.$$
   Thus the function ${\cal
   F}_{D}(K(iy,.)\phi)(z)$ is positive.\\ On the other hand  by taking $ y = x_j$
   we deduce from (6.14) the following relation
$$
   \sum_{n=0}^{\infty}||x_n||^{4p}{\cal
   F}_{D}(K(ix_j,.)\phi)(x_n) \geq ||x_j||^{4p}{\cal
   F}_{D}(K(ix_j,.)\phi)(x_j).
   \eqno{(6.16)}$$
But from the relation (6.15) we have $${\cal
   F}_{D}(K(ix_j,.)\phi)(x_j) = \int_{ \mathbb{R}^d}|K(ix_j,t|^2 \phi(t)\omega_k(t) dt.$$
By applying Lemma 6.1 there exist positive constants $C$ and $X$
such that for $||x_j|| \geq X$ we have $${\cal
   F}_{D}(K(ix_j,.)\phi)(x_j) \geq \frac{C}{
   ||x_j||^{2\gamma+d}}.$$
 From this inequality and (6.14),(6.16) we obtain for $||x_j|| \geq X$:
$$ \langle K(ix_j,.)\triangle_k^{2p}U,\phi\rangle \geq C
|x_j|^{4p-2\gamma-d}.$$
   Thus
$$\langle K(ix_j,.)\triangle_k^{2p}U,\phi\rangle \to +\infty, \;
as \; j \to +\infty.$$ This contradicts (6.14). Hence the
distribution $U$ is not given by a function $g\omega_k$ with $g$
in $ \mathcal{E}( \mathbb{R}^d).$\\ \noindent{\bf{Proposition
6.2.}} Let $S$ be in $ \mathcal{E'}( \mathbb{R}^d)$. If $S$ is
hypoelliptic then $S$ satisfies the ii) of the $H$-property.\\
\noindent{\bf{Proof}}\\ \hspace*{5mm} Suppose  that the $ii)$ of
the $H$-property does not hold. Then there exists a  sequence
$(z_n)_{n \in \mathbb{N}} \subset \mathbb{C}^d$ and a positive
constant $M$ such that for all $n \in \mathbb{N}$, ${\cal
F}_{D}(S)(z_n) = 0$ and $|Im z_n| \leq M Log |z_n|$.
\\
\hspace*{5mm} Let $\phi$ be in $\mathcal{D}(\mathbb{R}^d)$.
According to Theorem 4.1 i) there exists $a \in \mathbb{N}$ such
that for very $p \in \mathbb{N}$ we can find $C_p > 0$ for which
$$\forall \, z \, \in \mathbb{C}^d, \; |{\cal F}_{D}(\phi)(z)|
\leq C_p e^{a||Im z|| - p Log (1 + ||z||)}. $$ If we take $p \in
\mathbb{N}$ such that $ p
> Ma + 2$, we get $$||z_n||^2 |{\cal F}_{D}(\phi)(z_n)| \leq C_p. \eqno{(6.17)}$$
 Let
$(a_n)_{n \in \mathbb{N}}$ be a complex sequence such that the
series $\sum_{n = 0}^\infty |a_n|$ is convergent.\\ We consider
the sequence $\{ \mathcal{V}_q\}_{q \in \mathbb{N}}$ of
distributions in $ \mathcal{D'}(\mathbb{R}^d)$ given by
$$\mathcal{V}_q = \sum_{ n = 0}^{q} a_nT_{||z_n||^2
K(iz_n,.)\omega_k}.$$ For all $q,r \in \mathbb{N}$ with $q > r$,
we have
%\end{document}
$$\begin{array}{lll}\langle \mathcal{V}_q,\phi \rangle - \langle
\mathcal{V}_r,\phi\rangle &=& \langle\sum_{ n = r}^{q} a_n
T_{||z_n||^2 K(iz_n,.)\omega_k},\phi\rangle\\ &=& \sum_{ n =
r}^{q} a_n ||z_n||^2 {\cal F}_D(\phi)(-z_n). \end{array}$$ Thus
using (6.17) we obtain $$|\langle \mathcal{V}_q,\phi \rangle -
\langle \mathcal{V}_r,\phi\rangle| \leq C_p \sum_{ n = r}^{q}
|a_n| \to 0, \; as \, r \to +\infty. \eqno{(6.18)}$$ Then
$$\langle \mathcal{V}_q,\phi\rangle  \to L(\varphi), \; as \, q
\to +\infty.$$ We deduce that $L$ is a distribution $ \mathcal{V}$
in $ \mathcal{D'}(\mathbb{R}^d)$ and $\mathcal{V}_q$ converges to
$\mathcal{V}$ in $ \mathcal{D'}(\mathbb{R}^d)$ as $q$ tends to
infinity. Thus $$\mathcal{V} = \sum_{ n = 0}^{\infty}
a_nT_{||z_n||^2 K(iz_n,.)\omega_k},\eqno{(6.19)}$$ and from (6.18)
we deduce that $$|\langle \mathcal{V},\phi\rangle|\leq C_p \sum_{
n = 0}^{\infty} |a_n|. \eqno{(6.20)}$$

On the other hand by making a proof similar to those which has
given the relation (6.11) we obtain $$
 \langle S*_D \mathcal{V},\phi\rangle  =
 \langle \sum_{ n=0}^{\infty} a_n||z_n||^2
{\cal F}_D(S)(z_n) {\cal F}_D(\phi)(-z_n) = 0. $$ Thus $$S*_D
\mathcal{V} = 0 $$ As $S$ is hypoelliptic, we deduce that the
distribution $\mathcal{V}$ is given by a function $f\omega_k$ with
$f$
 in $ \mathcal{E}(\mathbb{R}^d)$. Then we have
$$ \mathcal{V} =  T_{f\omega_k}.\eqno{(6.21)}$$ From (2.8) we have
$$\forall \, n \in \mathbb{N}, \; K(iz_n,0) = 1.$$ Thus for all
closed ball of center $o$ and radius $r > 0$ we have  $$\forall \,
n \in \mathbb{N}, \; \sup_{y \in B(o,r)}|K(iz_n,y)| \geq
1.\eqno{(6.22)}$$ On the other hand using (6.20) and (6.21) we
obtain $$\sup_{y \in B(o,r)}|f(y)| \leq C_p \sum_{ n = 0}^{\infty}
|a_n|.$$ Thus $$\forall \, n \in \mathbb{N}, \; \sup_{y \in
B(o,r)}||z_n||^2|K(iz_n,y)| \leq C_p.\eqno{(6.23)}$$ From this
relation and (6.22) we deduce that
 $$\forall \, n \in \mathbb{N}, \; ||z_n| \leq C_p.\eqno{(6.24)}$$
which is a contradiction with our choice of the sequence $(z_n)_{n
\in \mathbb{N}}$. This completes the proof.\\
 \noindent{\bf{Proposition 6.3.}} Let $S$ be in ${\cal E'}(\mathbb{R}^d)$.
 If $S$ satisfies the $H$-property,
then there exists  a parametrix for $S$, that is, there exist
$\mathcal{V}$ in ${\cal E'}(\mathbb{R}^d)$ and $\psi$ in
$\mathcal{D}(\mathbb{R}^d)$ such that $\delta = S*_D \mathcal{V} +
T_{\psi \omega_k}$, where $\delta$ represents the Dirac
functional.

\noindent{\bf{Proof}}\\ \hspace*{5mm} Using (5.44)
 the $H$-property can also be written in the form \\ \hspace*{5mm}i)
  There exist $A,M > 0$
   such that $
  |{\cal
   F}( \,^{t}V_k (S))(x)| \geq ||x||^{-A}$ for all $||x|| \geq M.$
\\ \hspace*{8mm}ii) $\lim_{||z|| \to
\infty, \, z \in {\cal Z}} \frac{||Im z||}{Log ||z||} = \infty,$
where ${\cal Z} = \{z \in \mathbb{C}^d, \; {\cal
   F}( \,^{t}V_k (S)(z) =
0.\}$

We see that the $H$-property is true for the distribution
$\,^{t}V_k (S)$ of $ \mathcal{E'}( \mathbb{R}^d)$ in the case of
the classical Fourier transform ${\cal F}$ on $ \mathbb{R}^d$.
Then from [11] there exists a parametrix for $\,^{t}V_k (S)$, that
is, there exist $\mathcal{V}_o$ in  $ \mathcal{E'}( \mathbb{R}^d)$
and $\psi_0$ in $ \mathcal{D}( \mathbb{R}^d)$ such that $$\delta =
\; \,^{t}V_k (S)* \mathcal{V}_0 + T_{\psi_0}.\eqno{(6.25)}$$ As
the operator $\,^{t}V_k$ is a topological isomorphism from $
\mathcal{E'}( \mathbb{R}^d)$ onto itself, we deduce from (6.25)
and (3.17) that $$\delta = \,^{t}V_k(S)*
\,^{t}V_k(\,^{t}V_k^{-1}(\mathcal{V}_0)) +
\,^{t}V_k(\,^{t}V_k^{-1}(T_{\psi_0})).$$ Thus  $$\delta =\,
\,^{t}V_k(S)* \,\,^{t}V_k(\mathcal{V}) +
\,\,^{t}V_k(T_{\psi\omega_k}).\eqno{(6.26)} $$ with
$$\,^{t}V_k^{-1}(\mathcal{V}_0) = \mathcal{V}, \; and \;
\,^{t}V_k^{-1}(\psi_0) = \psi.$$ The distribution $\mathcal{V}$
and the function $\psi$ belong respectively to $ \mathcal{E'}
(\mathbb{R}^d)$ and $ \mathcal{D}( \mathbb{R}^d)$.\\ On the other
hand from Proposition 5.2 we have $$\,^{t}V_k(S*_D \mathcal{V})
=\, \,^{t}V_k(S)
* \,^{t}V_k(\mathcal{V}).$$ Thus the relation (6.26) can also be written in
the form $$\,^{t}V_k^{-1}(\delta) = S*_D \mathcal{V} +
T_{\psi\omega_k}.$$ But $$\,^{t}V_k^{-1}(\delta) = \delta.$$ Thus
$$ \delta = S*_D \mathcal{V} + T_{\psi \omega_k}.$$

\noindent{\bf{Theorem 6.1.}} We assume that the distribution $S$
in ${\cal E'}(\mathbb{R}^d)$ is such that ${\cal Z} = \{z \in
\mathbb{C}^d, \; {\cal F}_{D}(S)(z) = 0\}$ is infinite. The
following assertions are equivalent. \\ \hspace*{5mm}i) $S$ is
hypoelliptic.\\\hspace*{5mm}ii) $S$ satisfies the
$H$-properties.\\ \hspace*{5mm}iii) There exists a parametrix for
$S$, that is, there exist $\mathcal{V}$ in ${\cal
E'}(\mathbb{R}^d)$ and $\psi$ in $\mathcal{D}(\mathbb{R}^d)$ such
that $\delta = S*_D \mathcal{V} + T_{\psi \omega_k}$.
\\
\noindent{\bf{Proof}}\\ \hspace*{5mm}From Propositions 6.1 and 6.2
it suffices to prove that iii)
   $\Longrightarrow$ i). Assume that the distribution $\mathcal{U}$ is in $\mathcal{D'}(\mathbb{R}^d)$
    and that
   $S*_D \mathcal{U}$ is given by a function $f\omega_k$, with $f$ in $\mathcal{E}(\mathbb{R}^d)$.\\
    From iii) we have $$\delta = S *_D \mathcal{V} + T_{\psi \omega_k},$$ with $\mathcal{V}$
     in  ${\cal
   E'}(\mathbb{R}^d)$ and $\psi$ in $\mathcal{D}(\mathbb{R}^d)$.\\Thus
   $$\begin{array}{lll}
   \mathcal{U}  &=& \mathcal{U} *_D \delta \\
    &=& \mathcal{U}*_D( S *_D \mathcal{V} + T_{\psi \omega_k}).
     \end{array}
     $$ Using the commutativity and
     the associativity of the Dunkl convolution product of distributions in
     $\mathcal{E'}(\mathbb{R}^d)$, we obtain

$$\begin{array}{lll}
   \mathcal{U}  &=&  \mathcal{V} *_D( S *_D \mathcal{U}) + \mathcal{U}*_D T_{\psi \omega_k}\\
   &=& \mathcal{V}*T_{f \omega_k} + \mathcal{U}*_D T_{\psi \omega_k}.
     \end{array}
     $$
By applying (5.39) we obtain $$\begin{array}{lll}
   \mathcal{U}  &=&  T_{(\mathcal{V}*_D f)\omega_k} +
   T_{(\mathcal{U}*\psi)\omega_k}\\ &=& T_{(\mathcal{V}*_D f + \mathcal{U}*_D\psi)\omega_k}.
     \end{array}
     $$
     But from Theorem 5.1 the function $\mathcal{V}*_D f + \mathcal{U}*_D\psi$ belongs
     to $ \mathcal{E}(\mathbb{R}^d)$.\\ Thus $S$ is
     hypoelliptic.\\
\noindent{\bf{Example 6.1.}} \\ \hspace*{5mm} We suppose that $d
\geq 2$ and we consider the equation $$\triangle_k \mathcal{U} =
\mathcal{V},$$ $\mathcal{U}$ and $\mathcal{V}$ are distributions
in $ \mathcal{D'}( \mathbb{R}^d)$. We say that the Dunkl Laplacian
$\triangle_k$ is hypoelliptic if all solution $\mathcal{U}$ is
given by a function $f\omega_k$ with $f$ in $ \mathcal{E}(
\mathbb{R}^d)$ whenever $\mathcal{V}$ is given by a function
$g\omega_k$ with $g$ in $ \mathcal{E}( \mathbb{R}^d)$. As we have
$$\triangle_k \mathcal{U} = (\triangle_k \delta)*_D \mathcal{U},$$
where $\delta$ is the Dirac distribution on $ \mathbb{R}^d$.\\
Then the hypoellipticity of $ \triangle_k$ is equivalent to the
hypoellipticity of the distribution $ \triangle_k \delta$ in $
\mathcal{E'}( \mathbb{R}^d)$ given by
$$\langle\triangle_k\delta,\varphi\rangle =
\langle\delta,\triangle_k\varphi\rangle = \triangle_k \varphi(o),
\; \varphi \in  \mathcal{E}( \mathbb{R}^d).$$ The relation (2.8)
implies that $$\forall \, z \in \mathbb{C}^d, \; {\cal
F}_D(\triangle_k \delta)(z) = \sum_{j=1}^d z_j^2.\eqno{(6.27)}$$
i) From (6.27) we deduce that $$\forall \, x \in \mathbb{R}^d, \;
{\cal F}_D(\triangle_k \delta)(x) = ||x||^2.$$ Thus for $||x||
\geq 1$ we have  $$ |{\cal F}_D(\triangle_k \delta)(x)| \geq
||x||^{-1}.\eqno{(6.28)}$$ ii) The relation (6.27) implies also
that $$ \mathcal{Z} = \{z \in \mathbb{C}^d, \; {\cal
F}_D(\triangle_k \delta)(z) = 0\} = \{(x,y) \in \mathbb{R}^d\times
\mathbb{R}^d; \; ||x|| = ||y|| \; and \, \langle x,y\rangle =
0.\}$$ Thus $$\lim_{||z|| \to \infty, \, z \in {\cal Z}}
\frac{||Im z||}{Log ||z||} = \lim_{||y|| \to \infty, }
\frac{||y||}{Log (2^{\frac{1}{2}}||y||)} = +\infty.\eqno{(6.29)}
$$ The relations (6.28),(6.29) show that the distribution
$\triangle_k \delta$ satisfies the $H$-property. Thus Theorem 6.1
implies that the distribution $\triangle_k \delta$ is
hypoelliptic. \\ The Dunkl Laplacian is then hypoelliptic. This
result was first proved in \cite{MT1t} by another method.\\
% \end{document}

{ \Large{{\bf References}}}
\begin{enumerate}
\bibitem{B}{\bf Belhadji.M. and Betancor.J.J.(2002).}
{\em Hypoellipticity of Hankel Convolution equations on Schwartz
distibutions, Arch.Math. 79, 188-196. }

\bibitem{BBM}{\bf Betancor.J.
J., Betancor.J.D. and M\'endez.J.R.M. (2004) } {\em Hypoelliptic
Jacobi Convolution operators on Schwartz distibutions, Positivity
8,  407-422. }

\item {\bf Chazarain. J. and Piriou. A. (1982)}. Introduction to the
theory of linerar partial differential equations. North Holland
Publishing Company - Amsterdam, New-York. Oxford.
\item {\bf van Diejen. J. F. (1997)} . Confluent hypergeometric
orthogonal polynomials related to the rational quantum Calogero
system with harmonic confinement.Comm. Math. Phys., {\bf 188},
467-497.
\item {\bf Dunkl. C. F. (1989).} Differential-difference operators
associated to reflection groups. Trans. Amer.Math. Soc, {\bf 311},
167-183
\item {\bf Dunkl.C.F. (1991)}. Integral kernels with reflection group
invariance. Can. J. Math., {\bf 43}, 1213-1227.
\item {\bf Dunkl. C.F. (1992)}. Hankel transform associated to finite
reflection groups. Contemp. Math., {\bf 138}, 123-138.
\bibitem{E}{\bf Ehrenpreis.L. (1960).  }{\em Solution for
some problems of divisions IV. Invertible and elliptic operators,
Amer.J.Math. 82,  522-588 .}
\item {\bf Heckman. G.J. (1991)}. An elementary approach to the
hypergeometric shift operators of Opdam. Invent. Math., {\bf 103},
341-350.
\item {\bf Hikami. K. (1996).} Dunkl operators formalism for quantum
many-body problems associated with classical root systems. J.
Phys. Soc. Japan, {\bf 65}, 394-401.
\bibitem{Hor}{\bf H\"ormander.L. (1961). }{\em
Hypoelliptic convolution equations. Math.Scand. 9,  178-181. }
\item {\bf Humphreys. J.E. (1990)}. Reflection groups and Coxeter groups.
Cambridge Univ. Press., Cambridge, England.
\item {\bf de Jeu. M.F.E (1993)}. The Dunkl transform. Invent. Math.,
{\bf 113}, 147-162.
\item {\bf de Jeu. M.F.E. (1994)}.Dunkl operators. Thesis, University of
Amsterdam.
\item {\bf Kakei. S. (1996)}. Common algebraic structure for the
Calogero-Sutherland models. J. Phys. A, {\bf 29}, 619-624.
\item {\bf Lapointe. M. and Vinet, L. (1996)}. Exact operator solution
of the Calogero-Sutherland model. Comm. Math. Phys., {\bf 178},
425-452.
\bibitem{MT1t}
{\bf Mejjaoli.H. and Trim\`eche.K. (2001).} {\em On a mean value
property associated with
 the Dunkl Laplacian operator and applications.
    Integ. Transf. and Special Funct., Vol.12, $N^{0}$ 3,  279-301.}

\item {\bf R\"osler.M. and  Voit.L. (1998)}. Marcov processes related with
Dunkl operators. Adv. Appl. Math. 21, 575 -643
\item {\bf R\"osler. M. (1999)}. Positivity of Dunkl's intertwining
operator. Duke. Math. J., {\bf 98}, 445-463.
\item{\bf R\"osler. M. (2003)}. A positive radial product formulas for the Dunkl
kernel. Trans.Amer.Math.Soc.,{\bf 355}, 2413-2438.
\item{\bf Schwartz.L.(1966).} Th\'eorie des distributions. Hermann.
Editeurs des Sciences et des arts.
\item{\bf Thangavelu.S and Xu.Y. } Convolution operator and
maximal function for Dunkl transform. Arxiv:math.CA / 0403049.
\item {\bf Trim\`eche. K. (2001)}. The Dunkl intertwining operator on
spaces of functions and distributions and integral representation
of its dual. Integ. Transf. and Special Funct., {\bf 12} (4),
349-374.
\item {\bf  Trim\`eche.K. (2001)}. Generalized harmonic analysis and
wavelet packets. Gordon and Breach Science Publishers.
\item {\bf Trim\`eche.K. (2002)}. Paley-Wiener theorems for the Dunkl
transform and Dunkl translation operators. Integ. Transf. and
Special Funct., {\bf 13}, 17-38.
\item{\bf Trim\`eche.K.} Inversion formulas and geometrical form
of Paley-Wiener-Schwartz theorem associated with the Dunkl
operators.  Preprint.2005.
\item {\bf Xu. Y. (1997)}. Orthogonal polynomials for a family of
product weight function on the sphere. Can. J. Math. 49, 175-192.
\end{enumerate}
 \end{document}